\documentclass[12pt]{article}

\usepackage{amsmath}
\usepackage{amssymb}
\usepackage{amsthm}
\usepackage{amsfonts}
\usepackage{amscd}
\usepackage{graphicx}
\usepackage[T1]{fontenc}    
\usepackage{hyperref}
\usepackage[nottoc]{tocbibind}  
\usepackage[usenames,dvipsnames]{color}
\usepackage{array} 
\usepackage{enumitem} 
\usepackage{subfig} 
\usepackage[numbers,sort&compress]{natbib}
\usepackage{fancyvrb} 
\usepackage{geometry}
\usepackage{titling} 
\usepackage{tikz}
\usepackage[medium]{titlesec}
\usepackage{etoolbox}       
\usepackage{arydshln} 

\theoremstyle{plain}
\newtheorem{lemma}{Lemma}
\newtheorem{proposition}{Proposition}
\newtheorem{theorem}{Theorem}
\newtheorem{corollary}{Corollary}

\newtheorem{algorithm}{Algorithm}
\newtheorem{theoremext}{Theorem}

\theoremstyle{definition}
\newtheorem{definition}{Definition}
\newtheorem{remark}{Remark}
\newtheorem{example}{Example}

\hypersetup{colorlinks,urlcolor=blue!50!red,citecolor=black!45!green,linkcolor=blue}
\makeatletter\patchcmd{\ttlh@hang}{\parindent\z@}{\parindent\z@\leavevmode}{}{}\patchcmd{\ttlh@hang}{\noindent}{}{}{}\makeatother 
\titlespacing*{\section}{0pt}{1mm}{1mm}
\titlespacing*{\subsection}{0pt}{1mm}{1mm}
\titlespacing*{\paragraph}{0pt}{1mm}{1mm}
\newcommand{\myspace}{\setlength{\abovedisplayskip}{1mm}\setlength{\belowdisplayskip}{0mm}}
\newenvironment{Mlist}{\begin{itemize}[topsep=0pt,itemsep=0pt,leftmargin=7mm]}{\end{itemize}}

\newcommand{\PRP}[1]{Proposition~\ref{prp:#1}}
\newcommand{\LEM}[1]{Lemma~\ref{lem:#1}}
\newcommand{\THM}[1]{Theorem~\ref{thm:#1}}
\newcommand{\COR}[1]{Corollary~\ref{cor:#1}}
\newcommand{\DEF}[1]{Definition~\ref{def:#1}}
\newcommand{\RMK}[1]{Remark~\ref{rmk:#1}}
\newcommand{\EQN}[1]{(\ref{eqn:#1})}
\newcommand{\SEC}[1]{\textsection\ref{sec:#1}}
\newcommand{\FIG}[1]{Figure~\ref{fig:#1}}
\newcommand{\TAB}[1]{Table~\ref{tab:#1}}
\newcommand{\EXM}[1]{Example~\ref{exm:#1}}

\newcommand{\ALG}[1]{Algorithm~\ref{alg:#1}}
\newcommand{\END}{\hfill \ensuremath{\vartriangleleft}}
\newcommand{\fig}[3]{{\includegraphics[height=#1cm, width=#2cm]{img/#3}}}
\newcommand{\Iff}{if and only if }
\newcommand{\st}{{such that }}
\newcommand{\Wlog}{without loss of generality }
\newcommand{\resp}{respectively}
\newcommand{\wrt}{with respect to }
\newcommand{\AN}{\emptyset}
\newcommand{\Ai}{{A_1}}
\newcommand{\Aii}{{A_2}}
\newcommand{\Aiii}{{A_3}}
\newcommand{\uAi}{\underline{A_1}}
\newcommand{\uAii}{\underline{A_2}}
\newcommand{\uAiii}{\underline{A_3}}
\renewcommand{\P}{{\mathbb{P}}}
\renewcommand{\S}{{\mathbb{S}}}
\newcommand{\C}{{\mathbb{C}}}
\newcommand{\R}{{\mathbb{R}}}
\newcommand{\Z}{{\mathbb{Z}}}
\newcommand{\Q}{{\mathbb{Q}}}

\newcommand{\cW}{{\mathcal{W}}}

\newcommand{\cI}{{\mathcal{I}}}

\newcommand{\cLL}{{\mathcal{L}}}
\newcommand{\set}[2]{\{#1 ~|~ #2\}}
\newcommand{\df}[1]{\textit{#1}}
\newcommand{\aut}{\operatorname{Aut}}
\newcommand{\bas}[1]{\langle #1\rangle}
\newcommand{\p}{\varepsilon}
\renewcommand{\l}{\ell}
\renewcommand{\k}{k}
\newcommand{\Mi}{\mathfrak{i}}
\renewcommand{\c}{\colon}
\newcommand{\dto}{\dashrightarrow}

\newcommand{\sng}{\operatorname{sng}}
\newcommand{\piP}{\tilde{\pi}}
\newcommand{\g}{\overline{g}}

\begin{document}
\myspace

\begingroup  
\centering\Large Surfaces that are covered by two pencils of circles\par
\centering\large Niels Lubbes\par
\centering\large\today\par
\endgroup

\begin{abstract}
We list up to M\"obius equivalence all possible degrees and embedding dimensions of 
real surfaces that are covered by at least two pencils of circles,
together with the number of such pencils.
In addition, we classify incidences between the contained circles, complex lines and 
isolated singularities.
Such geometric characteristics are encoded in the N\'eron-Severi lattices of such
surfaces and is of potential interest to geometric modelers and architects.
As an application we confirm Blum's conjecture in higher dimensional space
and we address the Blaschke-Bol problem
by classifying surfaces that are covered by hexagonal webs of circles.
In particular, we find new examples of such webs that cannot be embedded in 3-dimensional space.
%
\end{abstract}

\begingroup
\def\addvspace#1{\vspace{-1mm}}
\tableofcontents
\endgroup

\section{Introduction}
\label{sec:intro}

\begin{figure}[!ht]
\centering
\begin{tabular}{@{}c@{\hspace{3mm}}c@{\hspace{3mm}}c@{\hspace{3mm}}c@{\hspace{3mm}}c@{}}
\fig{2.5}{2.5}{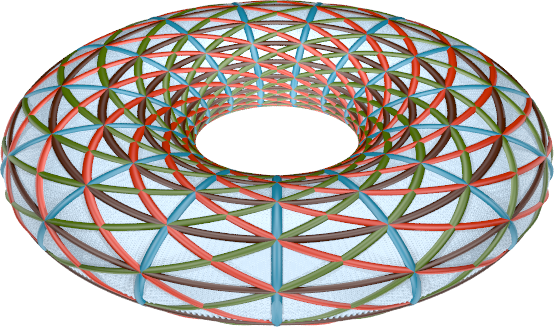}      & 
\fig{2.5}{2.5}{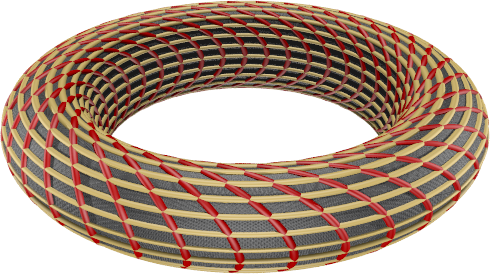} & 
\fig{2.5}{2.5}{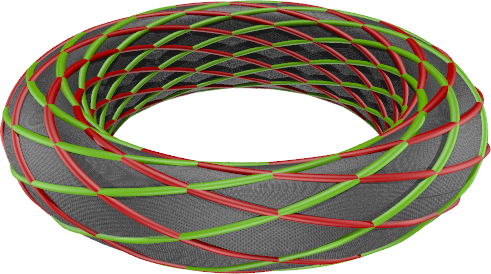} &
\fig{2.5}{2.5}{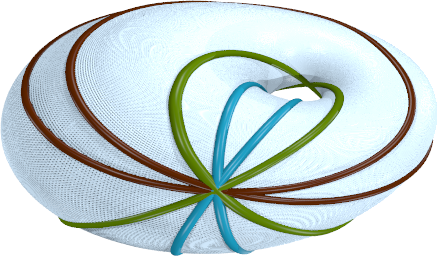}    &   
\fig{2.5}{2.5}{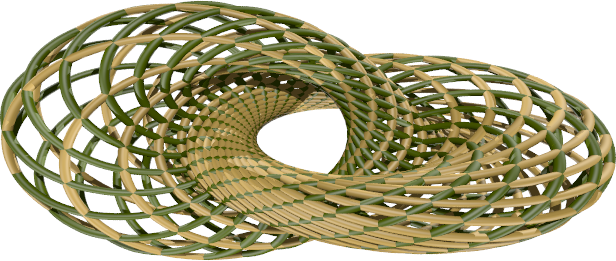}
\\
$(4,4,3)$ & $(5,4,3)$ &  & $(6,4,3)$ & $(2,8,3)$
\end{tabular}
\caption{celestial surfaces in $\R^3$ and their types.}
\label{fig:intro}
\end{figure}

Sir Christopher Wren discovered 
that a one-sheeted hyperboloid contains two lines through each point \citep[1669]{wrn1}
and he used his discovery for an {``engine designed for grinding hyperbolic lenses''} 
\citep[page~92]{bur1}.
We now know that a surface that is covered by two pencils of lines is 
either the plane, a hyperboloid of one sheet, or a hyperbolic paraboloid \citep[Lecture~16]{omni}.
In this article we consider
surfaces that are covered by two analytic pencils of circles instead of lines.
Such surfaces must be algebraic by \cite{conical}
(see also \citep[Theorem~4.1]{sko}).
With \df{surface} we shall therefore mean an irreducible algebraic surface.
We call a surface \df{$\lambda$-circled} 
if it contains no more and no less than
$\lambda$ circles through a general point and if the real points in the 
surface are not contained in a reducible curve.
For example, the leftmost ring cyclide in \FIG{intro} 
is a 4-circled quartic surface and appears in physical models such as 
Twistor theory~\citep[Fig. 33.14]{penrose}. 
A \df{celestial surface} is defined as a $\lambda$-circled surface \st $\lambda\geq 2$.
The name ``celestial'' is inspired by a model where planetary orbits 
are described by circles in 4-dimensional space \cite{bae}.
Celestial surfaces are of interest to
geometric modelers and architects~\cite{pot2,jut1, arch, zube, pet}.
In this section we introduce and motivate our main results, namely
\THM{T}, \THM{hex}, \THM{smooth}, \THM{B} and its corollaries
(see \SEC{results} for the proofs).

Let $\pi\c S^n\dto\R^n$ denote a stereographic projection  
from the $n$-dimensional unit-sphere $S^n\subset\R^{n+1}$.
Hipparchus of Nicaea (190--120 BCE) discovered that~$\pi$ 
sends circles to either circles or lines.
We say that a surface $Z\subset\R^n$ 
has \df{type} $(\lambda,d,n)$
if $\pi^{-1}(Z)\subset S^n$ is a $\lambda$-circled 
surface of degree~$d$ that 
is not contained in a hyperplane section (see \FIG{intro}).
Such types are invariant under M\"obius transformations of~$\R^n$.
We assume \Wlog that a celestial surface 
in~$\R^n$ is not contained in a hyperplane or hypersphere.
We know from \cite{conical} that if $\lambda\geq 2$, then $n<d\leq 8$
and either $\lambda\leq 10$ or $\lambda=\infty$.

\begin{theorem}
\label{thm:T}
The type $(\lambda,d,n)$ of a celestial surface, 
is one of the following and each listed type is realized by some surface:
$(2,8,n)$ for $3\leq n\leq 7$, $(3,6,5)$, $(3,6,4)$, $(2,6,5)$ $(2,6,4)$,
$(\lambda,4,3)$ for $2\leq \lambda \leq 6$, $(\infty,4,4)$ and $(\infty,2,2)$.
\end{theorem}

\begin{figure}[!ht]
\centering
\begin{tabular}{@{}c@{~~~}c@{~~~}c@{~~~}c@{}}
\fig{3}{3}{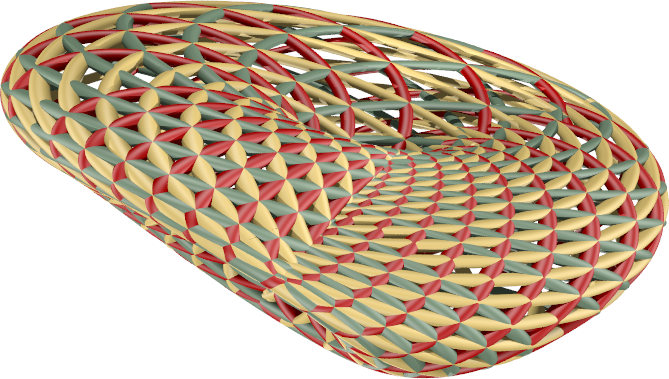}      &
\fig{3}{3}{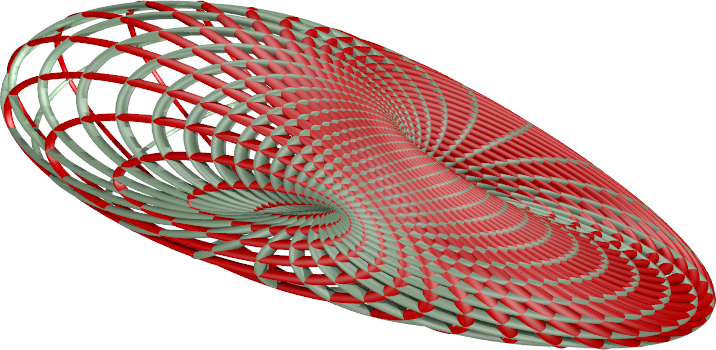}   & 
\fig{3}{3}{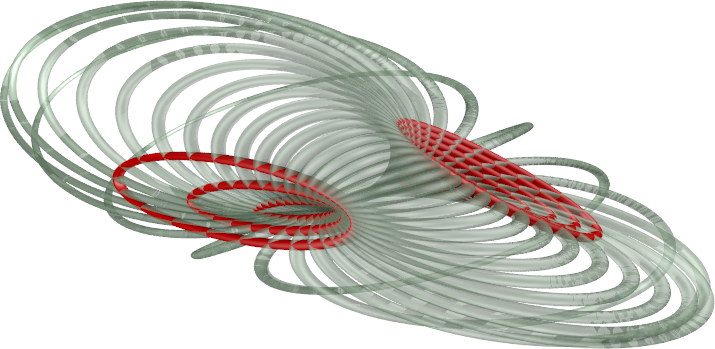}   &
\fig{3}{3}{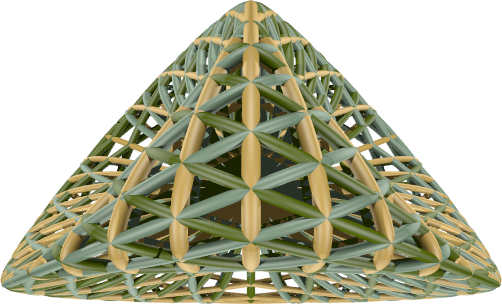}  
\\
$(3,6,5)$ & $(2,6,5)$ & $(2,6,5)$ & $(\infty,4,4)$
\end{tabular}
\caption{projected celestial surfaces are covered by ellipses.}
\label{fig:proj}
\end{figure}

In \FIG{proj} we see projections into $\R^3$ of celestial surfaces in $\R^n$ for $n>3$.

We recovered a result from \cite{kol2}, namely
that an $\infty$-circled surface is of type either $(\infty,4,4)$ or $(\infty,2,2)$.
It follows from \citep[Main~Theorem~1.1]{sko} that a celestial surface 
of type $(2,8,3)$
is M\"obius equivalent to either 
$\set{a+b\in \R^3}{a\in A,~b\in B}$ or
a stereographic projection
$\pi(\set{c \star d\in S^3}{a\in C,~b\in D})$, 
where $A,B\subset \R^3$ and $C,D\subset S^3$
are circles and $\star$ denotes the Hamiltonian product for the unit quaternions.

A surface of type $(\lambda,4,3)$ for some $\lambda>0$ is called a \df{Darboux cyclide}.
It follows from \RMK{dc} that this definition coincides with the definition in \citep[Section~2]{pot2}.
A systematic overview of equations for Darboux cyclides can be found in~\cite{tak2}.
The following corollary answers a question in \citep[Section~5]{pot2} 
and will be proven in the almost self-contained \SEC{type}.

\begin{corollary}
\label{cor:3}
A $\lambda$-circled surface in $\R^3$ \st $\lambda\geq 3$ 
is of degree at most four
and M\"obius equivalent to either a sphere or a Darboux cyclide.
\end{corollary}

The following corollary confirms the \df{Blum's conjecture},
which was known for manifolds in $\R^3$ 
that are homeomorphic to a torus~$S^1\times S^1$ \cite{tak1}
and for Darboux cyclides \citep[Remark~8]{pot2}.

\begin{corollary}
\label{cor:blum}
If a surface in $\R^n$ is $\lambda$-circled with $n\geq 2$,
then either $\lambda=\infty$ or $\lambda\leq 6$.
\end{corollary}

In order to refine and prove \THM{T} we 
propose to consider geometric aspects that happen at complex infinity.
To uncover this hidden structure we 
define a \df{real variety} $X$ to be a complex variety together with
an antiholomorphic involution $\sigma\c X\to X$ (see \citep[Section~I.1]{sil1} and \citep[Introduction]{moh})
and we denote its real points by $X\bas{\R}:=\set{p\in X}{\sigma(p)=p}$.
Such varieties can always be defined by polynomials with real coefficients \citep[Section~6.1]{serre}.

Points, curves, surfaces and projective spaces $\P^n$ are real algebraic varieties 
and maps between such varieties 
are compatible with the \df{real structure} $\sigma$ 
unless explicitly stated otherwise.
Moreover, a variety is irreducible
and a hypersurface or hyperplane section of a variety
inherits the real structure unless explicitly stated otherwise.
By default we assume that $X\bas{\R}$ for a surface $X$ is not contained in a 
reducible curve and that
the real structure $\sigma\c\P^n\to\P^n$ sends $x$ 
to $(\overline{x_0}:\ldots:\overline{x_n})$. 

As circles play a central role, it is natural to consider  
the \df{M\"obius quadric} for our space:
$
\S^n:=\set{x\in\P^{n+1}}{-x_0^2+x_1^2+\ldots+x_{n+1}^2=0}.
$
The \df{M\"obius transformations} of $\S^n$ are defined as biregular automorphisms $\aut(\S^n)$
and are linear so that $\aut(\S^n)\subset \aut(\P^{n+1})$.
If $V\subset \S^n$ is a variety, then we define $V(\R):=\gamma(V\bas{\R})$
where $\gamma\c \S^n\bas{\R}\to S^n$ is the isomorphism that sends $x$ to 
$\left(\frac{x_1}{x_0},\ldots,\frac{x_{n+1}}{x_0}\right)$.
Notice that $\pi^{-1}(\R^n)$
defines an isomorphic copy of $\R^n$ inside $\S^n(\R)$ \st the M\"obius transformations of $\S^n$
restrict via $\gamma$ and/or  $\pi$ 
to M\"obius transformations of $S^n$ and $\R^n$.

\begin{definition}
\label{def:model}
We call $C\subset \S^n$ a \df{circle} if $C(\R)\subset S^n$ is a circle
and a surface~$X\subset\S^n$ is \df{$\lambda$-circled} or \df{celestial}
if $X(\R)\subset S^n$ is as such.
We say that $X\subset\S^n$ is of \df{type} $(\lambda,d,n)$ if 
$X$ is a $\lambda$-circled surface of degree~$d$ that is not contained in a hyperplane section.
Notice that if $X\subset\S^n$ is a celestial surface, then its stereographic projection~$\pi(X(\R))$
is either a celestial surface in $\R^n$ or covered by lines. 
We call~$X$ in this case the \df{M\"obius model} of $\pi(X(\R))$.
We call $X\subset\S^3$ a \df{Darboux cyclide}
if it is the M\"obius model of a Darboux cyclide in~$\R^3$.
A \df{complex circle} is 
an irreducible complex conic $C\subset \S^n$.
\END
\end{definition}

An irreducible hypersurface $F\subset X\times \P^1$ is called
a \df{pencil} on a surface $X$ 
if $\pi_1(F)=X$ and $\pi_2(F)=\P^1$ 
for the projections of $F$ to its two factors.
We require that a \df{member} $F_i:=\pi_1(F\cap X\times\{i\})$ is a curve
for almost all points~$i\in \P^1$
and two pencils are equal if their members are the same.
The common complex points in the intersection $\cap_{i\in\P^1} F_i$ are called \df{base points}.
Notice that a point is real by convention, but that a ``base point'' is complex.
A \df{pencil of circles} is defined as a pencil whose general member is a circle.

The \df{circle graph} of a $\lambda$-circled surface $Z\subset\R^n$ \st $\lambda<\infty$ 
is defined as a labeled graph whose vertices correspond
to the pencils of circles that cover its M\"obius model $X\subset \S^n$. 
Each vertex is labeled with either $-$, $+$, $\times$ 
or no label,
if the pencil has two complex conjugate base points, one real base point, two real base points
and no base points, \resp.
Two vertices are connected by a solid or dashed labeled edge if 
general circles in the respective pencils intersect in two complex or real points.
We dash the edge if and only 
if at least one of the two intersection points coincides with a common base point
of the two pencils.  

\begin{example}
\label{exm:ring-graph}
Below we colored the vertices of the circle graph of the ring cyclide
to match the corresponding pencils of circles:
\begin{center}
\begin{tikzpicture}
\definecolor{wdred}{RGB}{153, 25, 25}
\definecolor{wlred}{RGB}{215, 91, 78}
\definecolor{wdblue}{RGB}{14 , 97,120}
\definecolor{wlblue}{RGB}{66 ,143,171}
\definecolor{wgreen}{RGB}{101,137, 23}
\definecolor{wblack}{RGB}{  0,  0,  0}
\definecolor{wbrown}{RGB}{ 86, 34,  0}
\node[inner sep=0pt] at (-3,0.5) {\fig{2.5}{3}{deg4-ring}}; 
\draw[draw=black] (0,1) -- (0, 0);
\draw[draw=black, fill=wlred!50 ] (0,1) circle [radius=3mm] node[yellow] {};
\draw[draw=black, fill=wgreen!50] (0,0) circle [radius=3mm] node[yellow] {};
\draw[draw=black, fill=wbrown!20] (1,1) circle [radius=3mm] node[black] {$-$};
\draw[draw=black, fill=wlblue!40] (1,0) circle [radius=3mm] node[black] {$-$};
\end{tikzpicture}
\end{center}
The two pencils of cospherical circles that cover the ring cyclide, 
are called \df{Villarceau circles} \citep[1848]{vil1}.
These circles can be found in a sculpture of a staircase in the Strasbourg cathedral,
which was built from 1176 until 1439 \citep[Fig. II.7.7]{ber2}.
The pencil of horizontal brown circles 
are parallel and thus all the circles in this pencil meet at complex conjugate base points at infinity.
The blue circles in the remaining pencil have complex conjugate base points at the axis of revolution.
\END
\end{example}

\begin{definition}
\label{def:web}
Suppose that~$\cW$ is a set of curves in a surface~$Z$.
We define $\cW_p:=\set{C\in\cW}{p\in C}$ for all~$p\in Z$
and let 
\[
V(\cW):=\set{v\in Z}{v \text{ is a point such that } |\cW_v|=3}.
\]
We call $\cW$ a \df{3-web} if $V(\cW)$ is not 
contained in some reducible curve.
We define $G(\cW)$ to be the graph with 
vertex set~$V(\cW)$ and 
labeled edge set
\[
\set{(\{v,w\},C)}{C\in\cW \text{ and } v,w\in C \text{ are pairwise distinct}}.
\]
We call a 3-web $\cW$ a \df{hexagonal web} if 
a general edge $\{p,q\}$ of the graph~$G(\cW)$
is contained in a subgraph as defined in \FIG{graph},
where the edge-labels $A,B,C,D,E,F,G,H,I\in\cW$ are pairwise distinct.
\END
\end{definition}

\begin{figure}[!ht]
\centering
\begin{tikzpicture}[scale=1.5]
\draw[black,          thick] (-1,0)--(-0.5,1) node [midway, left] {\footnotesize$H$} 
                                   --(0.5,1) node [midway, above] {\footnotesize$F$}
                                   --(1,0) node [midway, right] {\footnotesize$D$}
                                   --(0.5,-1) node [midway, right] {\footnotesize$E$}
                                   --(-0.5,-1) node [midway, below] {\footnotesize$G$}
                                   --(-1,0) node [midway, left] {\footnotesize$I$};
\draw[black!20!green, very thick] (-1,0) --(0,0) node[midway, above=-2pt] {\footnotesize$A$} --(1,0) node[midway, below=-2pt] {\footnotesize$A$};
\draw[red,            very thick] (0.5,1)--(0,0) node[midway, right=-2pt] {\footnotesize$C$} --(-0.5,-1) node[midway, left=-2pt] {\footnotesize$C$};
\draw[cyan,           very thick] (-0.5,1)--(0,0) node[midway, right=-2pt] {\footnotesize$B$} --(0.5,-1) node[midway, left=-2pt] {\footnotesize$B$}; 
\draw[black,fill=white] ( 0  , 0) circle [radius=0.2] node {$p$};
\draw[black,fill=white] (-1  , 0) circle [radius=0.2] node {$ $};
\draw[black,fill=white] ( 1  , 0) circle [radius=0.2] node {$q$};
\draw[black,fill=white] (-0.5, 1) circle [radius=0.2] node {$ $};
\draw[black,fill=white] ( 0.5, 1) circle [radius=0.2] node {$ $};
\draw[black,fill=white] (-0.5,-1) circle [radius=0.2] node {$ $};
\draw[black,fill=white] ( 0.5,-1) circle [radius=0.2] node {$ $};
\end{tikzpicture}
\caption{See \DEF{web}.}
\label{fig:graph}
\end{figure}

\begin{remark}
\label{rmk:web}
In order to understand \DEF{web} we suppose that $\cW$ is a hexagonal web on a surface
and that $\{p,q\}$ in \FIG{hex}a is a general edge of~$G(\cW)$.
In \FIG{hex}b we draw all the curves in~$\cW_p\cup\cW_q$
and we obtain at least two new intersection points $r$ and $s$.
In \FIG{hex}c we draw all the curves in $\cW_r\cup\cW_s$
and we obtain again at least two new intersection points. 
We repeat the last step one more time so that 
we obtain a closed hexagon as in \FIG{hex}d. 
\FIG{hex}e is an example for the case 
that $\cW$ is not a hexagonal web.  
We refer to \cite{nilov} for more information.
\END
\end{remark}

\begin{figure}[!ht]
\centering
\newcommand{\pp}[1] {\draw[draw=black, fill=gray!20, line width=0.1mm] #1 circle [radius=0.5mm];}

\def\tt{0.3}

\newcommand{\cL}[5] {\draw[cyan]           plot [smooth, tension=\tt] coordinates {#1 #2 #3 #4 #5};}
\newcommand{\dL}[4] {\draw[cyan]           plot [smooth, tension=\tt] coordinates {#1 #2 #3 #4};}
\newcommand{\eL}[3][] {\draw[cyan] #2 to [#1] #3;}
\newcommand{\cR}[5] {\draw[red]            plot [smooth, tension=\tt] coordinates {#1 #2 #3 #4 #5};}
\newcommand{\dR}[4] {\draw[red]            plot [smooth, tension=\tt] coordinates {#1 #2 #3 #4};}
\newcommand{\eR}[3][] {\draw[red] #2 to [#1] #3;}
\newcommand{\cH}[6] {\draw[black!20!green] plot [smooth, tension=\tt] coordinates {#1 #2 #3 #4 #5 #6};}
\newcommand{\dH}[4] {\draw[black!20!green] plot [smooth, tension=\tt] coordinates {#1 #2 #3 #4};}
\newcommand{\eH}[3][] {\draw[draw=black!20!green] #2 to [#1] #3;}

\newcommand{\cHH}[6] {\draw[black!20!green,densely dotted] plot [smooth, tension=\tt] coordinates {#1 #2 #3 #4 #5 #6};}

\def\d{0.05}

\def\Ax{-1              }\def\Ay{1              }
\def\Bx{-0.5            }\def\By{1              }
\def\Cx{0               }\def\Cy{1              }
\def\Dx{0.5-\d-\d       }\def\Dy{1              }
\def\Ex{1-\d-\d         }\def\Ey{1              }
\def\Fx{-1+\d+\d        }\def\Fy{-1             }
\def\Gx{-0.5+\d+\d      }\def\Gy{-1             }
\def\Hx{0               }\def\Hy{-1             }
\def\Ix{0.5-\d-\d-\d    }\def\Iy{-1             }
\def\Jx{1               }\def\Jy{-1             }
\def\Kx{-1              }\def\Ky{0.5+\d+\d      }
\def\Lx{-1              }\def\Ly{0.25+\d+\d     }
\def\Mx{-1              }\def\My{0+\d+\d        }
\def\Nx{-1              }\def\Ny{-0.25-\d-\d    }
\def\Ox{-1              }\def\Oy{-0.5           }
\def\Px{1               }\def\Py{0.5+\d+\d      }
\def\Qx{1               }\def\Qy{0.25-\d        }
\def\Rx{1               }\def\Ry{0              }
\def\Sx{1               }\def\Sy{-0.25          }
\def\Tx{1               }\def\Ty{-0.5           }
\def\Ux{-0.25-\d        }\def\Uy{0.25+\d+\d     }
\def\Vx{0.25-\d         }\def\Vy{0.25+\d        }
\def\Wx{0.5             }\def\Wy{0              }
\def\Xx{0.25            }\def\Xy{-0.25          }
\def\Yx{-0.25+\d        }\def\Yy{-0.25-\d       }
\def\Zx{-0.5-\d            }\def\Zy{0              }
\def\ZZx{\Zx}\def\ZZy{\Zy}
\def\OOx{0}\def\OOy{0}

\def\A{(\Ax,\Ay)}
\def\B{(\Bx,\By)}
\def\C{(\Cx,\Cy)}
\def\D{(\Dx,\Dy)}
\def\E{(\Ex,\Ey)}
\def\F{(\Fx,\Fy)}
\def\G{(\Gx,\Gy)}
\def\H{(\Hx,\Hy)}
\def\I{(\Ix,\Iy)}
\def\J{(\Jx,\Jy)}
\def\K{(\Kx,\Ky)}
\def\L{(\Lx,\Ly)}
\def\M{(\Mx,\My)}
\def\N{(\Nx,\Ny)}
\def\O{(\Ox,\Oy)}
\def\P{(\Px,\Py)}
\def\Q{(\Qx,\Qy)}
\def\R{(\Rx,\Ry)}
\def\S{(\Sx,\Sy)}
\def\T{(\Tx,\Ty)}
\def\U{(\Ux,\Uy)}
\def\V{(\Vx,\Vy)}
\def\W{(\Wx,\Wy)}
\def\X{(\Xx,\Xy)}
\def\Y{(\Yx,\Yy)}
\def\Z{(\Zx,\Zy)}
\def\OO{(\OOx,\OOy)}
\def\ZZ{(\ZZx,\ZZy)}

\def\CAP{(0,-1.2)}

\begin{tabular}{ccccc}
\begin{tikzpicture}
\cHH \R \W \OO \Z \ZZ \M;
\pp{\OO};\pp{\W}
\node at (\OOx,\OOy-0.3) {\footnotesize$p$};
\node at (\Wx,\Wy-0.3) {\footnotesize$q$};
\node[white] at \CAP {\tiny hexagonal web};
\end{tikzpicture}
&
\begin{tikzpicture}
\cL \A \U \OO \X \J;
\dL \B \V \W \T;
\cR \E \V \OO \Y \F;
\dR \P \W \X \G;
\cH \R \W \OO \Z \ZZ \M;
\draw[line width=0.7mm, draw=cyan] \W -- \V;
\draw[line width=0.7mm, draw=red]  \W -- \X;
\pp{\OO};\pp{\W};
\pp{\V};\pp{\X};
\node[white] at \CAP {\tiny hexagonal web};
\end{tikzpicture}
&
\begin{tikzpicture}
\cL \A \U \OO \X \J;
\dL \B \V \W \T;
\cR \E \V \OO \Y \F;
\dR \P \W \X \G;
\cH \R \W \OO \Z \ZZ \M;
\dH \Q \V \U \L;
\dH \S \X \Y \N;
\draw[line width=0.7mm, draw=cyan] \W -- \V;
\draw[line width=0.7mm, draw=red]  \W -- \X;
\draw[line width=0.7mm, draw=black!20!green] \V -- \U;
\draw[line width=0.7mm, draw=black!20!green] \X -- \Y;
\pp{\OO};\pp{\W};
\pp{\V};\pp{\X};
\pp{\U};\pp{\Y};
\end{tikzpicture}
&
\begin{tikzpicture}
\cL \A \U \OO \X \J;
\dL \B \V \W \T;
\dL \K \Z \Y \I;

\cR \E \V \OO \Y \F;
\dR \D \U \ZZ \O;
\dR \P \W \X \G;

\cH \R \W \OO \Z \ZZ \M;
\dH \Q \V \U \L;
\dH \S \X \Y \N;

\draw[line width=0.7mm, draw=cyan] \W -- \V;
\draw[line width=0.7mm, draw=red]  \W -- \X;
\draw[line width=0.7mm, draw=black!20!green] \V -- \U;
\draw[line width=0.7mm, draw=black!20!green] \X -- \Y;
\draw[line width=0.7mm, draw=red]  \U -- \ZZ;
\draw[line width=0.7mm, draw=cyan]  \Y -- \Z;

\pp{\OO};\pp{\W};
\pp{\V};\pp{\X};
\pp{\U};\pp{\Y};
\pp{\Z};\pp{\ZZ};

\end{tikzpicture}
&
\def\Ax{-1              }\def\Ay{1               }%
\def\Bx{-0.5            }\def\By{1               }%
\def\Cx{0               }\def\Cy{1               }%
\def\Dx{0.5-\d-\d       }\def\Dy{1               }%
\def\Ex{1-\d-\d-\d-\d   }\def\Ey{1               }%
\def\Fx{-1+\d+\d+\d+\d  }\def\Fy{-1              }%
\def\Gx{-0.5+\d+\d+\d+\d}\def\Gy{-1              }%
\def\Hx{0               }\def\Hy{-1              }%
\def\Ix{0.5-\d-\d-\d    }\def\Iy{-1              }%
\def\Jx{1               }\def\Jy{-1              }%
\def\Kx{-1              }\def\Ky{0.5+\d+\d       }%
\def\Lx{-1              }\def\Ly{0.25+\d+\d+\d+\d}%
\def\Mx{-1              }\def\My{0+\d+\d         }%
\def\Nx{-1              }\def\Ny{-0.25-\d-\d     }%
\def\Ox{-1              }\def\Oy{-0.5            }%
\def\Px{1               }\def\Py{0.5+\d+\d+\d+\d }%
\def\Qx{1               }\def\Qy{0.25-\d         }%
\def\Rx{1               }\def\Ry{0               }%
\def\Sx{1               }\def\Sy{-0.25           }%
\def\Tx{1               }\def\Ty{-0.5            }%
\def\Ux{-0.25-\d        }\def\Uy{0.25+\d+\d      }%
\def\Vx{0.25-\d         }\def\Vy{0.25+\d         }%
\def\Wx{0.5             }\def\Wy{0               }%
\def\Xx{0.25            }\def\Xy{-0.25           }%
\def\Yx{-0.25+\d        }\def\Yy{-0.25-\d        }%
\def\Zx{-0.5+\d+\d      }\def\Zy{0-\d            }%
\def\ZZx{-0.5-\d-\d-\d}\def\ZZy{0}%
\def\OOx{0}\def\OOy{0}%
\begin{tikzpicture}
\cL \A \U \OO \X \J;
\dL \B \V \W \T;
\dL \K \Z \Y \I;
\eL[out=-55, in=135]{\M}{(\Hx,\Hy)};
\eL[out=-45, in=135]{\C}{(\Rx,\Ry)};

\cR \E \V \OO \Y \F;
\dR \D \U \ZZ \O;
\dR \P \W \X \G;
\eR[out=200, in=45]{\R}{(\Hx,\Hy)};
\eR[out=220, in=45]{\C}{(\Mx,\My)};

\cH \R \W \OO \Z \ZZ \M;
\dH \Q \V \U \L;
\dH \S \X \Y \N;
\eH[out=190, in=15]{\P}{(\Kx,\Ky)};
\eH[out=190, in=-25]{\T}{(\Ox,\Oy)};

\draw[line width=0.7mm, draw=cyan] \W -- \V;
\draw[line width=0.7mm, draw=red]  \W -- \X;
\draw[line width=0.7mm, draw=black!20!green] \V -- \U;
\draw[line width=0.7mm, draw=black!20!green] \X -- \Y;
\draw[line width=0.7mm, draw=red]  \U -- \ZZ;
\draw[line width=0.7mm, draw=cyan]  \Y -- \Z;  

\pp{\OO};\pp{\W};
\pp{\V};\pp{\X};
\pp{\U};\pp{\Y};
\pp{\Z};\pp{\ZZ};
\end{tikzpicture}
\\
{\bf a} & {\bf b} & {\bf c} & {\bf d} & {\bf e}
\end{tabular}
\caption{See \RMK{web}.}
\label{fig:hex}
\end{figure}

\begin{remark}
\label{rmk:web2}
Discrete realizations of hexagonal webs lead to nice triangularizations of the underlying surface. 
We translate \citep[Theorem~18]{pot2} using the concept of circle graphs: 
{\it If $Z\subset\R^3$ is a Darboux cyclide, then 
three vertices of its circle graph form a hexagonal web \Iff
the vertices are not contained in the following labeled subgraph of the circle graph:}
\begin{tikzpicture}[scale=0.8] 
\draw[draw=black] ( 0, 0) -- ( 1, 0);
\draw[draw=black] ( 2, 0) -- ( 3, 0);
\draw[draw=black, fill=white] ( 0, 0) circle [radius=3mm] node[black] {};
\draw[draw=black, fill=white] ( 1, 0) circle [radius=3mm] node[black] {};
\draw[draw=black, fill=white] ( 2, 0) circle [radius=3mm] node[black] {};
\draw[draw=black, fill=white] ( 3, 0) circle [radius=3mm] node[black] {};
\end{tikzpicture}.
Notice that each vertex in the latter subgraph corresponds to a base point free pencil.
\END
\end{remark}

\begin{example}
\label{exm:ring-web}
Up to symmetries of the circle graph, there are two hexagonal webs of circles on a ring cyclide: 
\begin{center}
\begin{tabular}{c@{\hspace{2cm}}c}
\begin{tikzpicture}
\node[inner sep=0pt] at (-2.5,0.5) {\fig{2.5}{3}{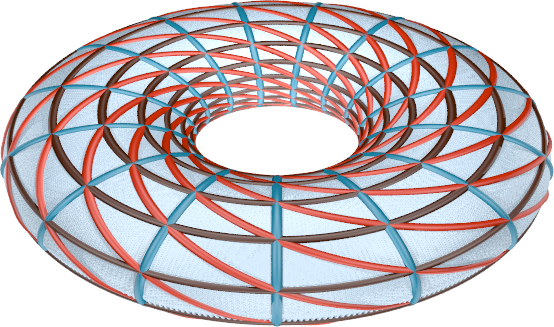}}; 
\definecolor{mred}{RGB}{215, 91, 78}
\definecolor{mblue}{RGB}{66 ,143,171}
\definecolor{mgreen}{RGB}{101,137, 23}
\definecolor{mbrown}{RGB}{ 86, 34,  0}
\draw[draw=black] (0,1) -- (0, 0);
\draw[draw=black, fill=mred!50 ] (0,1) circle [radius=3mm] node[yellow] {};
\draw[draw=black, fill=white] (0,0) circle [radius=3mm] node[yellow] {};
\draw[draw=black, fill=mbrown!40] (1,1) circle [radius=3mm] node[black] {$-$};
\draw[draw=black, fill=mblue!40] (1,0) circle [radius=3mm] node[black] {$-$};
\end{tikzpicture}
&
\begin{tikzpicture}
\node[inner sep=0pt] at (-2.5,0.5) {\fig{2.5}{3}{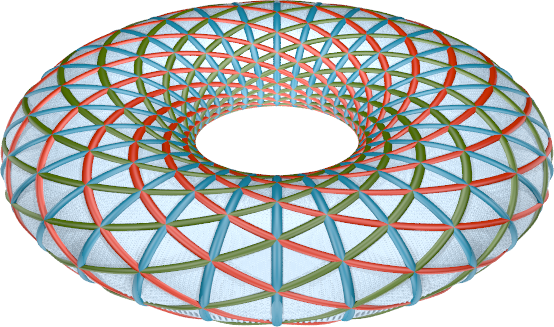}}; 
\definecolor{mred}{RGB}{215, 91, 78}
\definecolor{mblue}{RGB}{66 ,143,171}
\definecolor{mgreen}{RGB}{101,137, 23}
\definecolor{mbrown}{RGB}{ 86, 34,  0}
\draw[draw=black] (0,1) -- (0, 0);
\draw[draw=black, fill=mred!50 ] (0,1) circle [radius=3mm] node[yellow] {};
\draw[draw=black, fill=mgreen!50] (0,0) circle [radius=3mm] node[yellow] {};
\draw[draw=black, fill=white] (1,1) circle [radius=3mm] node[black] {$-$};
\draw[draw=black, fill=mblue!40] (1,0) circle [radius=3mm] node[black] {$-$};
\end{tikzpicture}
\end{tabular}
\end{center}
We colored the three vertices corresponding to the web.
\END
\end{example}

The following result addresses
the \df{Blaschke-Bol problem}, which refers to the 
classification problem for hexagonal webs of circles
\citep[\textsection 3, Aufgabe~1, page~31]{bla1}.

\begin{theorem}
\label{thm:hex}
If $Z\subset \R^n$ is a $\lambda$-circled surface  
\st $\lambda\geq 3$, then $Z$ is covered 
by a hexagonal web of circles.
%
\end{theorem}

See \FIG{proj} for linear projections 
of hexagonal webs of circles on celestial surfaces of types $(\infty,4,4)$ and $(3,6,5)$.

We will now introduce some concepts
from algebraic geometry such as ``smooth model'' and ``N\'eron-Severi lattice''.
Although these concepts are well-known to algebraic geometers, 
we would like to convince also non-experts of its usefulness. 
For example, we will see that the N\'eron-Severi lattice of 
a celestial surface encodes its circle graph. 

The \df{smooth model} of a surface $X\subset\P^n$ is a birational morphism $\varphi\c Y\to X$ 
from a nonsingular surface~$Y$, that does not contract complex $(-1)$-curves. 
See \citep[Theorem~2.16]{kolsing} for the existence and uniqueness of the smooth model. 

Suppose that $X\subset\S^n$ is a surface with smooth model~$Y\to X$
\st $Y$ is isomorphic 
to $\P^1\times\P^1$ blown up in either zero, two or four 
sufficiently general complex points that are left invariant as a set by the real structure $\sigma$.
With \df{sufficiently general} is meant 
that at most two centers of blowup are contained 
in a complex fiber of a projection of $\P^1\times\P^1$ to its first or second factor.
Under this assumption, $Y$ is an example of a weak del Pezzo surface (see \RMK{dp}).
Let $V$ denote the vector space of all forms of bidegree (2,2) on $\P^1\times\P^1$.
An \df{anticanonical model} of~$Y$ is the image of a birational map
$\P^1\times\P^1\dto X_N\subset\P^m$, 
whose components form a basis for the $(m+1)$-dimensional 
subspace of $V$ defined by the forms that 
vanish at the centers of blowup.
If $X$ is a degree preserving linear projection of~$X_N$, then
we show in \SEC{alg} how this point of view can be translated into 
an algorithm for constructing parametrizations 
and thus visualizations of celestial surfaces.

\begin{remark}
\label{rmk:S1S1}
If the real structure $\sigma\c \P^1\times\P^1\to \P^1\times\P^1$ maps $(x_0:x_1\,;\,x_2:x_3)$
to either 
$(\overline{x_0}:\overline{x_1}\,;\,\overline{x_2}:\overline{x_3})$
or
$(\overline{x_2}:\overline{x_3}\,;\,\overline{x_0}:\overline{x_1})$,
then $\P^1\times\P^1$ is isomorphic to $\S^1\times\S^1$ and~$\S^2$, \resp. 
\END
\end{remark}

We call $X\subset\S^n$ a \df{Veronese surface}
if it is isomorphic to the image of the biregular isomorphism $\P^2\to X_N\subset\P^5$,
whose components form a basis for the vector space of quadratic forms on~$\P^2$.

\begin{theorem}
\label{thm:smooth}
If $X\subset\S^n$ is celestial surface of type $(\lambda,d,n)$ with $n\geq 3$
and smooth model~$Y\to X$, 
then $d$, $\lambda$ and $Y$ are characterized by one of the rows in \TAB{smooth},
where the centers of blowup are in sufficiently general position.
If $\lambda<\infty$, then $X$ is either an anticanonical model of~$Y$
or a degree preserving linear projection of this anticanonical model.
If $\lambda=\infty$, then $X$ is a Veronese surface.
\end{theorem}

\begin{table}[!ht]
\caption{See \THM{smooth}.}
\label{tab:smooth}
\centering
\begin{tabular}{c@{\hspace{8mm}}c@{\hspace{8mm}}l}
$d$ & $\lambda$  & smooth model \\\hline
$8$ & $2$        & $\S^1\times\S^1$\\
$6$ & $\leq 3$   & $\S^1\times\S^1$ blown up in a pair of complex conjugate points \\
$4$ & $\infty$   & $\P^2$\\
$4$ & $\leq 6$   & $\S^1\times\S^1$ blown up in two pairs of complex conjugate points \\
$4$ & $2$        & $\S^2$ blown up in in two pairs of complex conjugate points \\
$4$ & $2$        & $X$ itself and $X(\R)$ consist of two disjoint spheres \\
\end{tabular}
\end{table}

We remark that if $X(\R)$ is not connected, then $X$ is not real birational
to the the plane $\P^2$ \citep[VI.6.5]{sil1}.

\begin{figure}[!ht]
\centering
\begin{tabular}{@{}c@{\hspace{5mm}}c@{\hspace{5mm}}c@{\hspace{5mm}}c@{}}
\fig{3}{3}{deg4-blum-2}           &
\fig{3}{3}{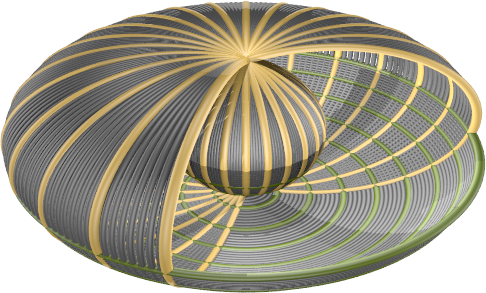}      & 
\fig{3}{3}{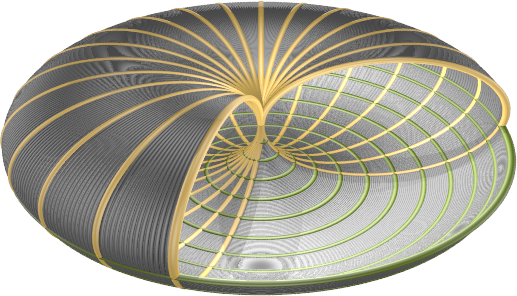}         &
\fig{3}{3}{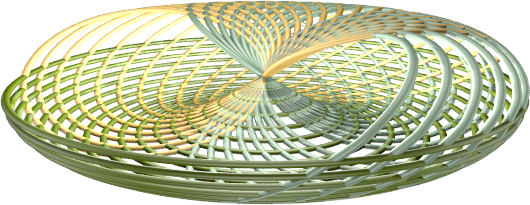}  
\end{tabular}
\caption{Blum cyclide, CO cyclide, CY cyclide and CH1 cyclide.}
\label{fig:quadric}
\end{figure}

\begin{definition}[names]
\label{def:names}
If $X\subset\S^3$ is of type $(4,4,3)$, $(5,4,3)$ or $(6,4,3)$ 
and $X(\R)$ is smooth, then $X$
is called a \df{ring cyclide}, \df{Perseus cyclide} and \df{Blum cyclide}, \resp~(see \FIG{intro}).
The latter two names have been introduced in \cite{blum1} and \cite{drei}, \resp.
We call $X\subset\S^3$
a \df{S1 cyclide} or \df{S2 cyclide} 
if $X(\R)$ is homeomorphic to $S^2$ and the disjoint union $S^2\,\dot{\cup}\, S^2$, \resp~(see \citep[Figure~12]{pot2}).
We use the following mnemonics for names of quadric surfaces 
whose equations are up to Euclidean similarity as in \TAB{quad}:

{\centering
\begin{tabular}{llll}
E = elliptic/ellipsoid      & P = parabolic/paraboloid   &  O = cone     \\ 
C = circular                & H = hyperbolic/hyperboloid &  Y = cylinder \\ 
\end{tabular}
}

For example, we call $X\subset\S^3$ a
\df{CH1 cyclide} if it is the M\"obius model of a Circular Hyperboloid of 1 sheet.
The \df{CO cyclide} and \df{CY cyclide} are also known as \df{spindle cyclide} 
and \df{horn cyclide}, \resp~(see \FIG{quadric}). 
We call a celestial of type $(2,8,n)$, $(3,6,n)$ or $(2,6,n)$ for some $n>0$,
\df{dS}, \df{dP6} and \df{wdP6}, \resp~(see \FIG{intro} and \FIG{proj}). 
If the M\"obius model of a surface $Z\subset \R^3$ is
a CH1 cyclide, then we call $Z$ also a \df{CH1 cyclide}.
Similarly, for the other names.
We remark that ``wdP'' stands for ``weak del Pezzo surface'' (see \RMK{dp})
and ``dS'' stands for ``double Segre surface'' 
\citep[8.4.1]{dol}.
\END
\end{definition}

Suppose that $\varphi\c Y\to X$ is the smooth model of a surface~$X\subset\P^n$.

The \df{N\'eron-Severi lattice} $N(X)$ is 
an additive group  defined by the divisor classes on $Y$ up to numerical equivalence.
This group comes with an unimodular intersection product $\cdot$
and a unimodular involution $\sigma_*\c N(X)\to N(X)$ induced by the real structure $\sigma\c X\to X$.
We denote by $\aut N(X)$ the group automorphisms that are compatible 
with both $\cdot$ and~$\sigma_*$.

The \df{class} $[C]\in N(X)$ of a complex curve $C\subset X$ 
is defined as the divisor class of 
the one-dimensional part of its complex preimage $\varphi^{-1}(C)$ minus 
the components that are contracted by the smooth model $\varphi\c Y\to X$.
We consider the following subsets of $N(X)$.
\begin{Mlist}
\item $B(X)$ denotes the set of classes of irreducible complex curves $C\subset Y$
\st $C$ is contracted by $\varphi$ to a complex point in $X$, and

\item $G(X)$ denotes the set of classes of complex irreducible conics in $X$
that are not components of the singular locus of~$X$.
\end{Mlist}
We call $W\subset B(X)$ a \df{component} if it defines 
a maximal connected subgraph of the graph with vertex set $B(X)$ and edge set $\set{(a,b)}{a\cdot b>0}$.
We write $c\cdot W\succ 0$ for $c\in N(X)$ and $W\subset N(X)$, if there exists $w\in W$ \st $c\cdot w>0$.
Similarly, we define $c\cdot W\nprec 0$, if there does not exists $w\in W$ \st $c\cdot w<0$.

\begin{proposition}
\label{prp:C}
Suppose that $X\subset\S^n$ is a celestial surface that is not $\infty$-circled.
\begin{Mlist}

\item[\bf a)] General circles $C,C'\subset X$ are members of the same pencil of circles
\Iff 
$[C]=[C']$, $[C]\in G(X)$ and $\sigma_*([C])=[C]$.

\item[\bf b)] General circles $C,C'\subset X$ intersect in two complex points 
\Iff 
there exists 
$2-[C]\cdot [C']$
components $W\subset B(X)$ \st $[C]\cdot W\succ0$ and $[C']\cdot W\succ0$. 

\item[\bf c)] The base points of a pencil of circles with member $C\subset X$ are
in one-to-one correspondence to the set of components 
$W\subset B(X)$ \st $[C]\cdot W\succ 0$. 
The base point is real \Iff $\sigma_*(W)=W$.



\end{Mlist}
\end{proposition}

We will consider N\'eron-Severi lattices that are generated by
$\bas{\l_0,\l_1,\p_1,\ldots,\p_r}_\Z$ for some $r\geq 0$
\st the nonzero intersections between 
its generators are $\l_0\cdot\l_1=1$ and $\p_1^2=\ldots=\p_r^2=-1$.
We define explicit coordinates for five different 
unimodular involutions $\sigma_*$ that act on such lattices (see below for $g_3$):
\begin{gather*}
\def\arraystretch{1.2}
\begin{array}{rl}
 A_0 :& r=0,~\sigma_*(\l_0)=\l_0,~ \sigma_*(\l_1)=\l_1,                                            \\
 A_1 :& r=2,~\sigma_*(\l_0)=\l_0,~ \sigma_*(\l_1)=\l_1,~ \sigma_*(\p_1)=\p_2,                      \\
2A_1 :& r=4,~\sigma_*(\l_0)=\l_0,~ \sigma_*(\l_1)=\l_1,~ \sigma_*(\p_1)=\p_2, \sigma_*(\p_3)=\p_4, \\
3A_1 :& r=4,~\sigma_*(\l_0)=\l_1,~ \sigma_*(\p_1)=\p_2,~ \sigma_*(\p_3)=\p_4,                      \\
D_4  :& r=4,~\sigma_*(\l_0)=g_3,~ \sigma_*(\l_1)=\l_1,~ \sigma_*(\p_i)=\l_1-\p_i \text{ for } 1\leq i\leq 4. 
\\
\end{array}
\end{gather*}

\begin{remark}
The names $A_1$, $2A_1$, $3A_1$ and $D_4$ correspond to the Dynkin types of 
root subsystems associated to $\sigma_*$ and these types are invariant under $\aut N(X)$.
See \cite{wal1} or the proof of \LEM{sigma} for details.
\END
\end{remark}

We use the following shorthand notation for elements in $B(X)$ and $G(X)$:
\[
\begin{array}{@{}l@{\hspace{1cm}}l@{\hspace{1cm}}l@{}}
b_1:=\p_1-\p_3, & b_{ij}:=\l_0-\p_i-\p_j   & b_0:=\l_0+\l_1-\p_1-\p_2-\p_3-\p_4, \\
b_2:=\p_2-\p_4, & b_{ij}':=\l_1-\p_i-\p_j. &
\end{array}
\]
\[
\begin{array}{@{}l@{\hspace{1cm}}l@{\hspace{1cm}}l@{}}
g_0:=\l_0, & g_2:=2\,\l_0+ \l_1-\p_1-\p_2-\p_3-\p_4, & g_{ij}:=\l_0+\l_1-\p_i-\p_j,\\
g_1:=\l_1, & g_3:= \l_0+2\,\l_1-\p_1-\p_2-\p_3-\p_4. & 
\end{array}
\]

\begin{example}
\label{exm:ring-real}
If $X\subset \S^3$ is a ring cyclide with smooth model~$Y\to X$, 
then $Y$ is isomorphic to $\S^1\times\S^1$ blown up
in two pairs of complex conjugate points and $\sigma_*$
is up to $\aut N(X)$ equal to $2A_1$. 
Moreover,
$B(X)=\{b_{13},b_{24},b_{14}',b_{23}'\}$ and $G(X)=\{g_0,g_1, g_{12}, g_{34}\}$.
Notice that the vertices of the circle graph in \EXM{ring-graph} correspond
to the elements of~$G(X)$. 
It follows from \PRP{C} that $g_0$ is the class of a circle 
in a pencil that has to two complex base points corresponding 
to the components $\{b_{14}'\}$ and $\{b_{23}'\}$ of $B(X)$.
The classes of the Villarceau circles are $g_{12}$ and $g_{34}$.
\END
\end{example}

If 
$c\in N(X)$
and
$\Psi\subset N(X)$, then we write $c\sim \Psi$
if $c$ is up to permutation of the generators $\p_1,\ldots,\p_4$ and up to switching generators $\l_0$ and $\l_1$, equal
to an element in $\Psi$. 
For example, $g_{34} \sim \{g_{12}\}$, $g_0\sim\{g_1\}$ and $g_2\sim\{g_3\}$.

\begin{theorem}
\label{thm:B}
If $X\subset\S^n$ is celestial surface of type $(\lambda,d,n)$ \st $\lambda<\infty$, 
then $G(X)=\set{c\in N(X)}{ c\sim\{g_0,g_2,g_{12}\},~ c\cdot B(X)\nprec 0 }$.
Moreover, $d$, $\lambda$, $\sigma_*$, $B(X)$ and the name of $X$ 
correspond up to $\aut N(X)$ to exactly one row in \TAB{B}.
\end{theorem}

\begin{table}[!ht]
\caption{See \THM{B}.}
\label{tab:B}
\centering
\begin{tabular}{@{}c@{~~}c@{~~}c@{~~}l@{}c|c@{~~}c@{~~}c@{~~}l@{}c@{}}
$d$ & $\lambda$ & $\sigma_*$ & $B(X)$                              & name    &       $d$ & $\lambda$ & $\sigma_*$ & $B(X)$                               & name       \\\hline
$8$ & $2$       & $A_0 $     & $\emptyset$                         & dS      &       $4$ & $3$       & $2A_1$     & $\{b_1,b_2,b_{12}\}$                 & EY         \\
$6$ & $3$       & $A_1 $     & $\emptyset$                         & dP6     &       $4$ & $2$       & $2A_1$     & $\{b_1,b_2,b_{12},b_{13}',b_{24}'\}$ & CY         \\
$6$ & $2$       & $A_1 $     & $\{b_{12}\}$                        & wdP6    &       $4$ & $3$       & $2A_1$     & $\{b_{12},b_{34}\}$                  & EO         \\
$4$ & $6$       & $2A_1$     & $\emptyset$                         & Blum    &       $4$ & $2$       & $2A_1$     & $\{b_{12},b_{34},b_{13}',b_{24}'\}$  & CO         \\
$4$ & $5$       & $2A_1$     & $\{b_1,b_2\}$                       & Perseus &       $4$ & $2$       & $3A_1$     & $\{b_0\}$                            & EE/EH2     \\
$4$ & $4$       & $2A_1$     & $\{b_{13},b_{24},b_{14}',b_{23}'\}$ & ring    &       $4$ & $2$       & $3A_1$     & $\{b_{13},b_{24}'\}$                 & EP         \\
$4$ & $4$       & $2A_1$     & $\{b_{12}\}$                        & EH1     &       $4$ & $2$       & $3A_1$     & $\emptyset$                          & S1         \\
$4$ & $3$       & $2A_1$     & $\{b_{13},b_{24},b_{12}'\}$         & CH1     &       $4$ & $2$       & $D_4$      & $\emptyset$                          & S2         \\
$4$ & $2$       & $2A_1$     & $\{b_{12},b_{34}'\}$                & HP      &           &           &            &                                      &            \\
\end{tabular}                                                                       
\end{table}

\COR{graph} below generalizes \citep[Theorem~20, page 296]{cool} to $n>3$:
{\it if the circle graph of $X\subset\S^n$ contains an edge, 
then $X$ must be a Darboux cyclide.}
Recall from \PRP{C} that the circle graph of~$X$
is encoded by $B(X)$ and $G(X)$.

\begin{corollary}
\label{cor:graph}
If a celestial surface is not $\infty$-circled, 
then its circle graph is in \TAB{graph}.
\end{corollary}

\begin{table}[!ht]
\caption{See \COR{graph}.}
\label{tab:graph}
\centering
\begin{tabular}{@{}c@{\hspace{6mm}}c@{\hspace{6mm}}c@{\hspace{6mm}}c@{\hspace{6mm}}c@{\hspace{6mm}}c@{\hspace{6mm}}c@{}}
\begin{tikzpicture}[scale=0.8] 
\draw[draw=black, fill=white] ( 0, 1) circle [radius=3mm] node[black] {};
\draw[draw=black, fill=white] ( 0, 0) circle [radius=3mm] node[black] {};
\end{tikzpicture}
&
\begin{tikzpicture}[scale=0.8] 
\draw[draw=black, fill=white] ( 0.5, 1) circle [radius=3mm] node[black] {};
\draw[draw=black, fill=white] ( 0, 0) circle [radius=3mm] node[black] {};
\draw[draw=black, fill=white] ( 1, 0) circle [radius=3mm] node[black] {};
\end{tikzpicture}
&
\begin{tikzpicture}[scale=0.8] 
\draw[draw=black, fill=white] ( 0, 1) circle [radius=3mm] node[black] {};
\draw[draw=black, fill=white] ( 0, 0) circle [radius=3mm] node[black] {$+$};
\end{tikzpicture}
&
\begin{tikzpicture}[scale=0.8] 
\draw[draw=black] (-0.7, 1) -- (-0.7, 0);
\draw[draw=black] ( 0.0, 1) -- ( 0.0, 0);
\draw[draw=black] ( 0.7, 1) -- ( 0.7, 0);
\draw[draw=black, fill=white] (-0.7, 1) circle [radius=3mm] node[black] {};
\draw[draw=black, fill=white] ( 0, 1) circle [radius=3mm] node[black] {};
\draw[draw=black, fill=white] ( 0.7, 1) circle [radius=3mm] node[black] {};
\draw[draw=black, fill=white] (-0.7, 0) circle [radius=3mm] node[black] {};
\draw[draw=black, fill=white] ( 0, 0) circle [radius=3mm] node[black] {};
\draw[draw=black, fill=white] ( 0.7, 0) circle [radius=3mm] node[black] {};
\end{tikzpicture}
&
\begin{tikzpicture}[scale=0.8] 
\draw[draw=black] (-1.2, 1) -- (-1.2, 0);
\draw[draw=black] ( 0.2, 1) -- ( 0.2, 0);
\draw[draw=black, fill=white] (-1.2, 1) circle [radius=3mm] node[black] {};
\draw[draw=black, fill=white] ( 0.2, 1) circle [radius=3mm] node[black] {};
\draw[draw=black, fill=white] (-1.2, 0) circle [radius=3mm] node[black] {};
\draw[draw=black, fill=white] ( 0.2, 0) circle [radius=3mm] node[black] {};
\draw[draw=black, fill=white] (-0.5, 0.5) circle [radius=3mm] node[black] {$-$};
\end{tikzpicture}
&
\begin{tikzpicture}[scale=0.8] 
\draw[draw=black] (-1, 1) -- (-1, 0);
\draw[draw=black, fill=white] (-1, 1) circle [radius=3mm] node[black] {};
\draw[draw=black, fill=white] (-1, 0) circle [radius=3mm] node[black] {};
\draw[draw=black, fill=white] ( 0, 1) circle [radius=3mm] node[black] {$-$};
\draw[draw=black, fill=white] ( 0, 0) circle [radius=3mm] node[black] {$-$};
\end{tikzpicture}
&
\begin{tikzpicture}[scale=0.8] 
\draw[draw=black] ( 0, 1) -- ( 0, 0);
\draw[draw=black, fill=white] ( 0, 1) circle [radius=3mm] node[black] {};
\draw[draw=black, fill=white] ( 0, 0) circle [radius=3mm] node[black] {};
\end{tikzpicture}
\\
dS & dP6 & wdP6 & Blum & Perseus & ring & S1/S2
\end{tabular}
\\[2mm]
\begin{tabular}{@{}c@{\hspace{5mm}}c@{\hspace{5mm}}c@{\hspace{5mm}}c@{\hspace{5mm}}c@{\hspace{5mm}}c@{\hspace{5mm}}c@{\hspace{2mm}}c@{}}
\begin{tikzpicture}[scale=0.8] 
\draw[line width=1, draw=black] (-1, 1) -- (0, 1);
\draw[line width=1, draw=black, densely dotted] ( -1, 0) -- ( 0, 0);
\draw[draw=black, fill=white] (-1, 1) circle [radius=3mm] node[black] {};
\draw[draw=black, fill=white] ( 0, 1) circle [radius=3mm] node[black] {};
\draw[draw=black, fill=white] (-1, 0) circle [radius=3mm] node[black] {$+$};
\draw[draw=black, fill=white] ( 0, 0) circle [radius=3mm] node[black] {$+$};
\end{tikzpicture}
&
\begin{tikzpicture}[scale=0.8] 
\draw[line width=1, draw=black, densely dotted] ( -1, 0) -- ( 0, 0);
\draw[draw=black, fill=white] (-0.5, 1) circle [radius=3mm] node[black] {$-$};
\draw[draw=black, fill=white] (-1, 0) circle [radius=3mm] node[black] {$+$};
\draw[draw=black, fill=white] ( 0, 0) circle [radius=3mm] node[black] {$+$};
\end{tikzpicture}
&
\begin{tikzpicture}[scale=0.8] 
\draw[line width=1, draw=black, densely dotted] ( -1, 0) -- ( 0, 0);
\draw[draw=black, fill=white] (-1, 0) circle [radius=3mm] node[black] {$+$};
\draw[draw=black, fill=white] ( 0, 0) circle [radius=3mm] node[black] {$+$};
\end{tikzpicture}
&
\begin{tikzpicture}[scale=0.8] 
\draw[line width=1, draw=black] ( -1, 0) -- ( 0, 0);
\draw[draw=black, fill=white] (-0.5, 1) circle [radius=3mm] node[black] {$+$};
\draw[draw=black, fill=white] (-1, 0) circle [radius=3mm] node[black] {};
\draw[draw=black, fill=white] ( 0, 0) circle [radius=3mm] node[black] {};
\end{tikzpicture}
&
\begin{tikzpicture}[scale=0.8] 
\draw[draw=black, fill=white] ( 0, 1) circle [radius=3mm] node[black] {$+$};
\draw[draw=black, fill=white] ( 0, 0) circle [radius=3mm] node[black] {$-$};
\end{tikzpicture}
&
\begin{tikzpicture}[scale=0.8] 
\draw[line width=1, draw=black] ( -1, 0) -- ( 0, 0);
\draw[draw=black, fill=white] (-0.5, 1) circle [radius=3mm] node[black] {$\times$};
\draw[draw=black, fill=white] (-1, 0) circle [radius=3mm] node[black] {};
\draw[draw=black, fill=white] ( 0, 0) circle [radius=3mm] node[black] {};
\end{tikzpicture}
&
\begin{tikzpicture}[scale=0.8] 
\draw[draw=black, fill=white] (-1, 1) circle [radius=3mm] node[black] {$\times$};
\draw[draw=black, fill=white] (-1, 0) circle [radius=3mm] node[black] {$-$};
\end{tikzpicture}
&
\begin{tikzpicture}[scale=0.8] 
\draw[draw=black] ( 0, 1) -- ( 0, 0);
\draw[draw=black, fill=white] ( 0, 1) circle [radius=3mm] node[black] {};
\draw[draw=black, fill=white] ( 0, 0) circle [radius=3mm] node[black] {};
\end{tikzpicture}
\\
EH1 & CH1 & HP & EY & CY & EO & CO & EE/EH2/EP
\end{tabular}
\end{table}

Notice that if $X\subset\S^n$ is covered by a pencil of circles with a 
real base point, then the stereographic projection $\pi(X(\R))$ from this base point
is a surface in~$\R^n$ that is covered by lines.
\COR{line} below together with
\THM{smooth} and \SEC{alg}
addresses \citep[Problem~5.6]{sko}.
See also \citep[Theorem~1.1]{nil} for $n\leq 3$.

\begin{corollary}
\label{cor:line}
A surface of degree $\delta$ in $\R^n$ for $n\geq 3$ that
is not contained in a hyperplane section, 
and that contains $c\geq 1$ circles and $\ell\geq 1$ lines through a general point,
is characterized by a row in \TAB{line}.
\end{corollary}

\begin{table}[!ht]
\caption{See \COR{line}.}
\label{tab:line}
\centering
\begin{tabular}{c@{\qquad}c@{\qquad}c@{\qquad}c@{\qquad}l@{\qquad}l}
$n$ & $\delta$ & $c$ & $\ell$   & type           &  name          \\\hline
3   & 2      & $2$       & $2$      & $(4,4,3)$      &  EH1           \\
3   & 2      & $1$       & $2$      & $(3,4,3)$      &  CH1           \\
3   & 2      & $2$       & $1$      & $(3,4,3)$      &  EO or EY      \\
3   & 2      & $1$       & $1$      & $(2,4,3)$      &  CO or CY      \\
4   & 3      & $\infty$  & $1$      & $(\infty,4,4)$ &                \\
4   & 4      & $1$       & $1$      & $(2,6,4)$      &                \\
5   & 4      & $1$       & $1$      & $(2,6,5)$      &                \\
\end{tabular}
\end{table}

If $X\subset \S^n$ is a celestial surface,
then, by \LEM{BGE}, $N(X)$ encodes beside the circle graph
also partial information about complex lines in $X$ and the \df{singular locus} $\sng X$.
In addition to $B(X)$ and $G(X)$ we let $E(X)$ denote the classes of complex lines in $X$
and we use the following notation for its elements:
\[
e_i=\p_i,\qquad e_{ij}=\l_i-\p_j,\qquad e_i'=b_0+\p_i.
\]
A component $W\subset B(X)$ in \THM{B} defines a Dynkin graph of type 
either $\Ai$, $\Aii$ or $\Aiii$ and corresponds by \LEM{B} to  
an isolated double point that is a node, cusp or tacnode, \resp. 
We underline if the isolated singularity is real (for example $\uAi$)
and take formal sums to denote the disjoint union of singularities (for example $2\Ai+\uAi$). 

\begin{corollary}
\label{cor:sng}
If $X\subset\S^n$ is celestial surface of type $(\lambda,d,n)$
\st $n=d-1$, then 
$E(X)=\set{c\in N(X)}{c\sim\{e_1,e_{01},e_1'\},~ c\cdot B(X)\nprec 0 }$
and 
$\lambda$, $\sigma_*$, $\sng X$, $|E(X)|$
and the name of~$X$ 
correspond up to $\aut N(X)$ to a row in \TAB{sng}. 
\end{corollary}

\begin{table}[!ht]
\caption{See \COR{sng}.}
\label{tab:sng}
\centering
\begin{tabular}{@{}c@{\hspace{3mm}}cl@{\hspace{0mm}}cl|c@{\hspace{3mm}}cl@{\hspace{0mm}}cl@{}}
$\lambda$ & $\sigma_*$ & $\sng X$     & $|E(X)|$ & name    & $\lambda$ & $\sigma_*$ & $\sng X$      & $|E(X)|$  & name       \\\hline
$2$       & $A_0 $     & $\AN$        & $0$      & dS      & $3$       & $2A_1$     & $\uAiii$      & $4$       & EY         \\
$3$       & $A_1 $     & $\AN$        & $6$      & dP6     & $2$       & $2A_1$     & $\uAiii+2\Ai$ & $2$       & CY         \\
$2$       & $A_1 $     & $\uAi$       & $3$      & wdP6    & $3$       & $2A_1$     & $2\uAi$       & $8$       & EO         \\
$6$       & $2A_1$     & $\AN$        & $16$     & Blum    & $2$       & $2A_1$     & $2\uAi+2\Ai$  & $4$       & CO         \\
$5$       & $2A_1$     & $2\Ai$       & $8$      & Perseus & $2$       & $3A_1$     & $\uAi$        & $12$      & EE/EH2     \\
$4$       & $2A_1$     & $4\Ai$       & $4$      & ring    & $2$       & $3A_1$     & $\uAii$       & $8$       & EP         \\
$4$       & $2A_1$     & $\uAi$       & $12$     & EH1     & $2$       & $3A_1$     & $\AN$         & $16$      & S1         \\
$3$       & $2A_1$     & $\uAi+2\Ai$  & $6$      & CH1     & $2$       & $D_4$      & $\AN$         & $16$      & S2         \\
$2$       & $2A_1$     & $\uAii$      & $8$      & HP      &           &            &               &           &            \\
\end{tabular}                                                                       
\end{table}

\begin{example}
\label{exm:ring-line}
We continue with \EXM{ring-real}, where $X\subset \S^3$ 
is a ring cyclide and $B(X)=\{b_{13},b_{24},b_{14}',b_{23}'\}$
so that $E(X)=\{e_1,e_2,e_3,e_4\}$. 
The complex line with class $e_1$ intersects the complex line with class $e_3$ at the 
complex isolated singularity corresponding to the component $\{b_{13}\}$,
since $e_1\cdot b_{13}, e_3\cdot b_{12}>0$.
We find that $X$ contains two pairs of complex conjugate lines that intersect in two pairs of complex 
conjugate nodes. 
\END
\end{example}

For convenience of the reader we computed $E(X)$ and $G(X)$
from $B(X)$ in \THM{B} for all celestial Darboux cyclides.

\begin{corollary}
\label{cor:BEG}
If $X\subset\S^3$ is a celestial Darboux cyclide,
then its name together with $B(X)$, $E(X)$ and $G(X)$ are up to $\aut N(X)$
defined by a row in \TAB{BEG}.
\end{corollary}

\begin{table}[!ht]
\caption{See \COR{BEG}. 
A class is send by the unimodular involution $\sigma_*$ to 
itself if underlined and to
its left or right neighbor in the listing otherwise.
The dashed row dividers indicate that $\sigma_*$ is $2A_1$, $3A_1$ and $D_4$, \resp.}
\label{tab:BEG}
\centering
$
\begin{array}{@{}ll@{}} 
\text{name}      & B(X), E(X), G(X)
\\\hline
\text{Blum}      & \{\},                                                       \{e_1,e_2,e_3,e_4,e_{01},e_{02},e_{03},e_{04},e_{11},e_{12},e_{13},e_{14},e_4',e_3',e_2',e_1'\},    
\\               & \{\underline{g_0},\underline{g_1},\underline{g_{12}},\underline{g_{34}},\underline{g_2},\underline{g_3},g_{13},g_{24},g_{14},g_{23}\}
\\\text{Perseus} & \{b_1,b_2\},                                                \{e_3,e_4,e_{01},e_{02},e_{11},e_{12},e_4',e_3'\},                                                  \{\underline{g_0},\underline{g_1},\underline{g_{12}},\underline{g_2},\underline{g_3},g_{13},g_{24}\}
\\\text{ring}    & \{b_{13},b_{24},b_{14}',b_{23}'\},                          \{e_1,e_2,e_3,e_4\},                                                                                \{\underline{g_0},\underline{g_1},\underline{g_{12}},\underline{g_{34}}\}
\\\text{EH1}     & \{\underline{b_{12}}\},                                     \{e_1,e_2,e_3,e_4,e_{03},e_{04},e_{11},e_{12},e_{13},e_{14},e_2',e_1'\},                            
\\               & \{\underline{g_0},\underline{g_1},\underline{g_{34}},\underline{g_3},g_{13},g_{24},g_{14},g_{23}\}
\\\text{CH1}     & \{b_{13},b_{24},\underline{b_{12}'}\},                      \{e_1,e_2,e_3,e_4,e_{13},e_{14}\},                                                                  \{\underline{g_0},\underline{g_1},\underline{g_{34}},g_{14},g_{23}\}
\\\text{HP}      & \{\underline{b_{12}},\underline{b_{34}'}\},                 \{e_1,e_2,e_3,e_4,e_{03},e_{04},e_{11},e_{12}\},                                                    \{\underline{g_0},\underline{g_1},g_{13},g_{24},g_{14},g_{23}\}
\\\text{EY}      & \{\underline{b_{12}},b_1,b_2\},                             \{e_3,e_4,e_{11},e_{12}\},                                                                          \{\underline{g_0},\underline{g_1},\underline{g_3},g_{13},g_{24}\}
\\\text{CY}      & \{\underline{b_{12}},b_1,b_2,b_{13}',b_{24}'\},             \{e_3,e_4\},                                                                                        \{\underline{g_0},\underline{g_1}\}
\\\text{EO}      & \{\underline{b_{12}},\underline{b_{34}}\},                  \{e_1,e_2,e_3,e_4,e_{11},e_{12},e_{13},e_{14}\},                                                    \{\underline{g_0},\underline{g_1},\underline{g_3},g_{13},g_{24},g_{14},g_{23}\}
\\\text{CO}      & \{\underline{b_{12}},\underline{b_{34}},b_{13}',b_{24}'\},  \{e_1,e_2,e_3,e_4\},                                                                                \{\underline{g_0},\underline{g_1},g_{14},g_{23}\}
\\\hdashline
\text{EE/EH2}    & \{\underline{b_0}\},                                        \{e_1,e_2,e_3,e_4,e_{01},e_{12},e_{02},e_{11},e_{03},e_{14},e_{04},e_{13}\},                        
\\               & \{g_0,g_1,\underline{g_{12}},\underline{g_{34}},g_{13},g_{24},g_{14},g_{23}\}
\\\text{EP}      & \{b_{13},b_{24}'\},                                         \{e_1,e_2,e_3,e_4,e_{02},e_{11},e_{04},e_{13}\},                                                    \{g_0,g_1,\underline{g_{12}},\underline{g_{34}},g_{14},g_{23}\}
\\\text{S1}      & \{\},                                                       \{e_1,e_2,e_3,e_4,e_{01},e_{12},e_{02},e_{11},e_{03},e_{14},e_{04},e_{13},e_4',e_3',e_1',e_2'\},    
\\               & \{g_0,g_1,\underline{g_{12}},\underline{g_{34}},g_{13},g_{24},g_{14},g_{23},g_2,g_3\}
\\\hdashline
\text{S2}        & \{\},                                                       \{e_1,e_{11},e_2,e_{12},e_3,e_{13},e_4,e_{14},e_{01},e_1',e_{02},e_2',e_{03},e_3',e_{04},e_4'\},    
\\               & \{g_0,g_3,g_{12},g_{34},g_{13},g_{24},g_{14},g_{23},\underline{g_1},\underline{g_2}\}
\end{array}
$
\end{table}

\section{Smooth models}
\label{sec:high}

In \PRP{dp} we will characterize the possible
types $(\lambda,d,n)$ and smooth models
of a celestial surface $X\subset\S^n$.
For this purpose, 
we start by collecting known results from the 
literature and we prove \PRP{C}.

Let $X\subset\P^n$ be a surface with smooth model~$\varphi\c Y\to X$.
Its \df{linear normalization} $X_N\subset\P^m$ 
is defined as the image of~$Y$ via the map $\varphi_{h}$ associated to 
the linear equivalence class $h$ of the pullback to~$Y$ of a hyperplane section of~$X$.

\newpage
\begin{remark}
\label{rmk:XN}
The associated map $\varphi_{h}\c Y\to X_N\subset \P^m$ 
is compatible with the real structures
and $X_N$ is unique up to $\aut(\P^m)$ 
as a direct consequence of the definitions
(see \citep[Remark~II.7.8.1]{har} and \citep[Section~1.1.B]{laz}).
We have that $m\geq n$ and there exists a degree-preserving linear projection $\eta\c\P^m\to\P^n$ \st $\eta(X_N)=X$
(see \citep[Theorem~6]{conical}).
\END
\end{remark}

\begin{definition}
\label{def:blowup}
A surface $Y_{r+1}$ is the \df{blowup of $Y_1$ in $r$ complex points}
if there exists a birational morphism $\tau\c Y_{r+1}\to Y_1$
and complex blowups 
$\tau_i\c Y_{i+1}\to Y_{i}$ of 
complex points $p_i\in Y_i$
\st $\tau=\tau_1\circ\cdots\circ\tau_r$ (see \citep[Example~I.4.9.1]{har}).
We refer to $p_i$ as a \df{center of blowup}.
If $p_{i+1}$ lies on the complex $(-1)$-curve that is contracted by 
$\tau_i$, then we say that $p_{i+1}$ is \df{infinitely near} to~$p_i$.
If $C\subset Y_1$ is a complex curve, then its \df{strict transform} $C_i\subset Y_i$
for some~$i$
is defined as the Zariski closure of the preimage $(\tau_1\circ\cdots\circ\tau_{i-1})^{-1}(C)$
minus the complex curve components that are contracted by $\tau_1\circ\cdots\circ\tau_{i-1}$.
\END
\end{definition}

\begin{remark}
\label{rmk:dp}
Suppose that $X\subset\P^n$ is a surface with smooth model~$Y\to X$
\st $Y$ is complex isomorphic to 
the blown up of $\P^1\times\P^1$
in $0\leq r\leq 7$ sufficiently general points.
Notice that if $r>0$, then $Y$ is 
complex isomorphic to a blowup of $\P^2$ in $r+1$ complex points.
Let $\k\in N(X)$ denote the canonical class of~$Y$ 
(see \citep[Example~V.1.4.4]{har}).
If $-\k$ 
is the class of a hyperplane section of~$X$, 
then $Y$ is
a \df{weak del Pezzo surface} and we call 
its linear normalization~$X_N$ an \df{anticanonical model} of~$Y$ (see \citep[Definition~8.1.18]{dol}).
\END
\end{remark}

\begin{theoremext}[Schicho, 2001]
\label{thm:mc}
Suppose that $X\subset\P^n$ is a surface of degree~$d$
that contains 
$\lambda\geq 2$ complex conics
and no complex line
through a general complex point.
Let $Y\to X$ be its smooth model and let $\k$ denote its canonical class.
\begin{Mlist}
\item
If $\lambda<\infty$,
then $3\leq d\leq 8$,
$Y$ is complex isomorphic to 
$\P^1\times\P^1$ blown up in $8-d$ 
sufficiently general points,
$X_N\subset\P^d$
and
the class of a hyperplane section is equal to $-\k$.

\item
If $\lambda=\infty$, then $d=4$, 
$Y$ is isomorphic to $\P^2$,
$X_N\subset\P^5$
and the class of a hyperplane section is equal to $-\frac{2}{3}\k$.
\end{Mlist}
\end{theoremext}

\begin{proof}
Direct consequence of \citep[Theorems~5--8 and Proposition~1]{conical}
and \RMK{dp}. 
Notice that $-\frac{1}{3}\k$ is the class of a line in $\P^2$ \citep[Example~II.8.20.3]{har}.
\end{proof}

\begin{lemma}
\label{lem:BGE}
Suppose that $X\subset\P^m$ is a surface with smooth model~$Y\to X$
\st $Y$ is complex isomorphic to the blowup of $\P^1\times\P^1$ 
in~$r\in\{0,2,4\}$ 
sufficiently general complex points.
Let $\pi_1,\pi_2\c\P^1\times\P^1\to\P^1$ denote the first and second projection.
\begin{Mlist}
\item[\bf a)]
$N(X)$
is generated by $\bas{\l_0,\l_1,\p_1,\ldots,\p_r}_\Z$
\st
\begin{itemize}[topsep=0mm,itemsep=0mm,leftmargin=0mm]
\item[$\bullet$]
$\l_0$ and $\l_1$ are the pullbacks of the classes of 
the general fibers of $\pi_1$ and $\pi_2$, \resp,
\item[$\bullet$]
$\p_1,\ldots,\p_r$ are 
the pullbacks of 
the classes of complex $(-1)$-curves 
that contract to the centers $p_1,\ldots,p_r$ of blowup,
and
\item[$\bullet$]
the canonical class $\k$ of $X$ is equal to 
$-2\,\l_0-2\,\l_1+\p_1+\ldots+\p_r$.
\end{itemize}

\item[\bf b)]
$b\in B(X)$ \Iff there exists $1\leq i<j\leq 4$ \st either
\begin{itemize}[topsep=0mm,itemsep=0mm,leftmargin=0mm]
\item[$\bullet$] $b=b_{ij}$ so that $p_i$ and $p_j$ lie on strict transforms of 
a complex fiber of $\pi_1$,
\item[$\bullet$] $b=b_{ij}'$ so that  $p_i$ and $p_j$ lie on strict transforms of 
a complex fiber of $\pi_2$,
\item[$\bullet$] $b=b_0$ so that $p_1,\ldots,p_4$ lie on strict transforms
of a complex curve of bidegree~$(1,1)$, or
\item[$\bullet$] $b=\p_i-\p_j$ so that  $p_j$ is infinitely near to $p_i$.
\end{itemize}

\item[\bf c)]
$G(X)=\set{c\in N(X)}{c\sim\{g_0,g_2,g_{12}\},~ c\cdot B(X)\nprec 0}$. 

\item[\bf d)]
$E(X_N)=\set{c\in N(X)}{c\sim\{e_1,e_1',e_{01}\},~ c\cdot B(X)\nprec 0}\subset E(X)$.

\end{Mlist}
\end{lemma}

\begin{proof}
Assertion a) follows from \citep[Proposition~V.2.3 and Proposition~V.3.3]{har}.
The classes in $B(X)$ and $E(X)$ are 
classes of complex $(-2)$-curves and complex $(-1)$-curves, \resp~(see \citep[8.2.7 and 8.2.6]{dol}).
Assertions b), c) and d) follow from \citep[Proposition~1, Lemma~1 and Proposition~2]{web}
after converting to a different set of generators for $N(X)$. 
See alternatively \citep[Proposition~8.2.7 and Lemma~8.2.22]{dol}
for assertions b) and d).
For d), notice that a complex non-linear curve in the linear normalization $X_N$
may be linearly projected to a complex line in
the singular locus of~$X$.
\end{proof}

\begin{proof}[Proof of \PRP{C}.]
Let $\varphi\c Y\to X$ be the smooth model of a celestial surface $X\subset\S^n$
that is not $\infty$-circled
and suppose that $C,C'\subset X$ are general circles.
Notice that $[C]\in G(X)$
and that $\sigma_*([C])=[C]$ as a direct consequence of the definitions.
It follows from \THM{mc}, \LEM{BGE}
and the
Riemann-Roch and Kawamata-Viehweg vanishing theorems
that a pencil of complex circles that covers $X$ are the images of 
complex curves in~$Y$ that form a
one-dimensional complete linear series without base points.
In particular, $[C]=[C']$ \Iff $C$ and $C'$ are members of the same pencil,
and thus assertion~a) holds.
Now let $D,D'\subset Y$ denote the preimages of $C$ and $C'$ in~$Y$, \resp.
We find that $C$ and $C'$ 
intersect in at least $[C]\cdot[C']=|D\cap D'|$ complex ``moving intersection'' points. 
Suppose that there exists
complex points
$p\in D$ and $q\in D'$
\st 
$p,q\notin D\cap D'$
and
$\varphi(p)=\varphi(q)$.
Thus $C$ and $C'$
will intersect in addition to the moving intersection points
at the complex point $\varphi(p)$.
As $C$ and $C'$ are chosen generally,
we find that $p$ and $q$
must lie on a 
possibly reducible complex curve $R\subset Y$ \st each
of its irreducible complex
components is contracted to~$\varphi(p)$.
Therefore, $\varphi(p)$ must be a base point for each pencil
that has the circle $C$ or $C'$ as member.
We have proven assertion~b)
as $|C\cap C'|$ minus the moving intersection points is 
equal to the number of base points that are met by both $C$ and $C'$.
Notice that 
$R$ corresponds to a component $W\subset B(X)$
and thus $\sigma(R)=R$ \Iff $\sigma_*(W)=W$ \Iff
the base point $\varphi(p)$ is real.
This concludes the proof for assertion~c).
\end{proof}

\begin{lemma}
\label{lem:even}
If $X\subset \S^n$ is a celestial surface,
then its degree is even and $X(\R)$ does not contain lines.
\end{lemma}

\begin{proof}
Suppose by contradiction that $X$ is of odd degree. 
The intersection of $X\bas{\R}\subset \P^{n+1}$ with the hyperplane at infinity 
is an odd degree curve and thus non-empty. 
We arrived at a contradiction as $X(\R)\subset S^n$.
In particular, $X(\R)$ is compact and thus does not contain lines.
\end{proof}

An \df{(2/3)-anticanonical model of $\P^2$}
is defined as the image of a map associated to 
minus two-thirds times the canonical class of $\P^2$.

\begin{proposition}
\label{prp:dp}
If $X\subset\S^n$ is a celestial surface of type $(\lambda,d,n)$ for $n>2$
with linear normalization $X_N$,
then up to $\aut N(X)$ either
\begin{Mlist}
\item 
$d=8$, 
$\lambda=2$, 
$n\leq 7$, 
$X_N\subset\P^8$ is an anticanonical model of $\S^1\times\S^1$, 
$\sigma_*=A_0$,
$B(X)=\emptyset$ and
$G(X)=\{g_0,g_1\}$,
\item 
$d=6$, 
$\lambda\leq 3$, 
$n\leq 5$, 
$X_N\subset\P^6$ is an anticanonical model of 
$\S^1\times\S^1$ blown up in a pair of complex conjugate points,
$\sigma_*=A_1$,
$B(X)\subseteq\{b_{12}\}$ and 
$G(X)\subseteq\{g_0,g_1,g_{12}\}$,

\item 
$d=4$, 
$\lambda=\infty$, 
$n\leq 4$, the Veronese surface $X_N\subset\P^5$
is a $(2/3)$-anticanonical model of $\P^2$
and
$\sigma_*$ is the identity, or

\item 
$d=4$, 
$\lambda\leq 10$, $n=3$
and $X\subset \P^4$ is complex isomorphic 
to an anticanonical model 
of~$\P^1\times\P^1$ blown up in four sufficiently general points.
\end{Mlist}
\end{proposition}

\begin{proof}
If $\lambda=\infty$, then the assertion follows from \THM{mc}.
We henceforth assume that $\lambda<\infty$
and let $\varphi\c Y\to X$ be the smooth model of~$X$.
It follows from \THM{mc} and \LEM{even} that $d\in\{4,6,8\}$
so that the assumption of \LEM{BGE} holds.
Moreover, $X_N\subset\P^d$ is an anticanonical model of~$Y$.
Notice that if $g\in G(X)$ is the 
class of a circle then $\sigma_*(g)=g$.
For the following three if-statements we apply \LEM{BGE} and \PRP{C}.
If $d=4$, then $\lambda\leq 10$
and thus the proof is concluded by \THM{mc}.
If $d=8$, then $\set{g\in G(X)}{\sigma_*(g)=g}=\{g_0,g_1\}$
and thus $Y\cong \S^1\times\S^1$ by \RMK{S1S1} as asserted.
If $d=6$, then $2\leq \lambda\leq 3$
and we may assume \Wlog that $g_0$ and $g_1$ are the classes of circles.
Since $g_0\cdot g_1=1$, 
the fiber product of the associated maps $\varphi_{g_0},\varphi_{g_1}\c Y\to \P^1$ 
defines a birational morphism $\varphi_{g_0}\times \varphi_{g_1}\c Y\to \S^1\times\S^1$.
Hence, $Y$ must be isomorphic to $\S^1\times\S^1$ blown up in $8-d$ centers.
If $E\subset Y$ is a resulting complex $(-1)$-curve, then $-k\cdot [E]=1$ and 
thus $\varphi(E)$ is a complex line in~$X$, which must be non-real by \LEM{even}
so that the centers are complex conjugate as is asserted. 
\end{proof}

\begin{lemma}
\label{lem:B}
Suppose that 
$\varphi\c Y\to X_N$ defines the smooth model
of a linear normal surface $X_N$, \st 
$Y$ is complex isomorphic to the blowup of $\P^1\times\P^1$ in $0\leq r\leq 4$ 
sufficiently general points.
If $W\subset B(X)$ is a component, then
the union of curves $C\subset Y$ \st $[C]\in W$ 
is contracted by $\varphi$ to an isolated double point of $X_N$.
If the graph with vertex set $W$ 
and edge set $\set{(a,b)}{a\cdot b>0}$
has Dynkin type $\Ai$, $\Aii$ or $\Aiii$, then this 
double point is a node, cusp or tacnode, \resp.
\end{lemma}

\begin{proof}
See \citep[Proposition~8.1.10]{dol} and \citep[Theorem~8.2.28]{dol}.
\end{proof}

\section{Constructing celestial surfaces}
\label{sec:alg}

In this section we propose a method 
for constructing examples of celestial surfaces of given type. 
Each celestial surface can be constructed with this method
and we show that if $(d,n)\notin\{(6,3),(4,3)\}$,
then each type in \THM{T} is realized by a celestial surface.
However, our method leaves room for improvement and this will be formulated
as an open problem.

Each entry in \TAB{cert} is marked with a type~$(\lambda,d,n)$
and consists of three items ({\tt u},{\tt R},{\tt Q})
which encode a 
rational map $\mu\c \S^1\times\S^1\dto \P^d$,
a linear projection $\rho\c\P^d\dto\P^{n+1}$ and
a hyperquadric $Q\subset\P^{n+1}$. 
We will show that if $\nu\c\P^{n+1}\to\P^{n+1}$ 
is a projective transformation \st $\nu(Q)=\S^n$,
then $\nu\circ\rho\circ\mu$ 
is a birational maps of bidegree~$(2,2)$
that parametrizes a celestial surface of type~$(\lambda,d,n)$.
We start by explaining the encoding with an example (recall \RMK{S1S1}).

\begin{table}[!ht]
\caption{
See \LEM{uRQ}.
The list {\tt u}, matrix {\tt R} and symmetric matrix {\tt Q},
defines the map $\mu$, projection $\rho$ and hyperquadric $Q$, \resp.}
\label{tab:cert}
\centering
\begin{Verbatim}[fontsize=\tiny,baselinestretch=1,commandchars=\\\{\}]
\textcolor{red}{# type (3,6,5)}
u = [x^2*v^2-y^2*w^2,x^2*v*w+y^2*v*w,x^2*w^2+y^2*w^2,x*y*v^2-y^2*v*w,x*y*v*w-y^2*w^2,y^2*v*w+x*y*w^2,y^2*v^2+y^2*w^2]
R = [(1,0,0,0,0,0,0),(0,1,0,0,0,0,0),(0,0,1,0,0,0,0),(0,0,0,1,0,0,0),(0,0,0,0,1,0,0),(0,0,0,0,0,1,0),(0,0,0,0,0,0,1)] 
Q = [(0,0,2,0,5/2,3,4),(0,-4,0,-5/2,-1/2,-2,-7/2),(2,0,0,-5/2,2,0,15/2),(0,-5/2,-5/2,-8,1/2,-9/2,0),
     (5/2,-1/2,2,1/2,-6,1/2,-1),(3,-2,0,-9/2,1/2,-7,7/2),(4,-7/2,15/2,0,-1,7/2,0)]

\textcolor{red}{# type (3,6,4)}
u = [x^2*v^2-y^2*w^2,x^2*v*w+y^2*v*w,x^2*w^2+y^2*w^2,x*y*v^2-y^2*v*w,x*y*v*w-y^2*w^2,y^2*v*w+x*y*w^2,y^2*v^2+y^2*w^2]
R = [(0,0,0,1,0,1,1),(1,1,0,1,0,0,0),(0,1,1,0,1,0,1),(1,1,0,0,1,0,0),(1,1,1,1,1,1,0),(1,0,0,1,0,1,1)]
Q = [(-20,-9,-3/2,16,-1/2,-31/2),(-9,-2,7/2,-9,11/2,25/2),(-3/2,7/2,0,-14,19/2,29/2),(16,-9,-14,16,5,-21),
     (-1/2,11/2,19/2,5,-19,-3),(-31/2,25/2,29/2,-21,-3,25)]

\textcolor{blue}{# type (2,6,5)}
u = [ x^2*v^2 + y^2*v^2,x^2*v*w,x^2*w^2,x*y*v*w,x*y*w^2,y^2*v*w,y^2*w^2 ]
R = [(1,0,0,0,0,0,0),(0,1,0,0,0,0,0),(0,0,1,0,0,0,0),(0,0,0,1,0,0,0),(0,0,0,0,1,0,0),(0,0,0,0,0,1,0),(0,0,0,0,0,0,1)]
Q = [(0,0,4,0,5/2,0,7/2),(0,-8,0,-5/2,-1,-7/2,-4),(4,0,0,1,0,4,4),(0,-5/2,1,-8,0,-5/2,3),(5/2,-1,0,0,-8,-3,0),
     (0,-7/2,4,-5/2,-3,-7,0),(7/2,-4,4,3,0,0,0)] 

\textcolor{blue}{# type (2,6,4)}
u = [ x^2*v^2 + y^2*v^2,x^2*v*w,x^2*w^2,x*y*v*w,x*y*w^2,y^2*v*w,y^2*w^2 ]
R = [(1,0,0,1,1,1,0),(1,0,0,1,0,1,1),(0,1,1,1,0,1,0),(0,0,0,0,1,0,0),(0,1,0,1,1,1,0),(0,0,0,1,1,1,0)]
Q = [(-6,3,-7/2,11/2,14,-15/2),(3,0,19/2,19/2,-20,15/2),(-7/2,19/2,0,-1,0,9/2),(11/2,19/2,-1,-27,-14,33/2),
     (14,-20,0,-14,-12,3/2),(-15/2,15/2,9/2,33/2,3/2,-6)] 

\textcolor{magenta}{# type (2,8,7)}
u = [x^2*v^2,x^2*v*w,x^2*w^2,x*y*v^2,x*y*v*w,x*y*w^2,y^2*v^2,y^2*v*w,y^2*w^2]
R = [(1,0,0,0,0,0,0,0,0),(0,1,0,0,0,0,0,0,0),(0,0,1,0,0,0,0,0,0),(0,0,0,1,0,0,0,0,0),(0,0,0,0,1,0,0,0,0),
     (0,0,0,0,0,1,0,0,0),(0,0,0,0,0,0,1,0,0),(0,0,0,0,0,0,0,1,0),(0,0,0,0,0,0,0,0,1)] 
Q = [(0,0,9/2,0,-1/2,-7/2,9/2,-7/2,15/2),(0,-9,0,1/2,0,1/2,0,-7/2,4),(9/2,0,0,7/2,-1/2,0,3,-3,4),
     (0,1/2,7/2,-9,7/2,-7/2,0,-5/2,2),(-1/2,0,-1/2,7/2,-7,-1,5/2,-5,4),(-7/2,1/2,0,-7/2,-1,-8,3,-4,0),
     (9/2,0,3,0,5/2,3,0,0,4),(-7/2,-7/2,-3,-5/2,-5,-4,0,-8,0),(15/2,4,4,2,4,0,4,0,0)] 

\textcolor{magenta}{# type (2,8,6)}     
u = [x^2*v^2,x^2*v*w,x^2*w^2,x*y*v^2,x*y*v*w,x*y*w^2,y^2*v^2,y^2*v*w,y^2*w^2] 
R = [(0,0,0,1,1,1,0,0,0),(0,1,0,0,0,0,0,0,0),(1,0,1,1,0,0,0,1,0),(1,0,0,1,0,0,1,0,0),(0,0,0,0,0,0,1,0,1),
     (0,1,0,0,1,0,1,0,1),(0,0,0,1,1,0,1,0,1),(1,0,0,1,0,0,1,1,0)]
Q = [(-8,6,0,-4,3/2,-17/2,7,-9/2),(6,-20,3/2,1,-3/2,11,-9/2,0),(0,3/2,0,7,3/2,-3/2,4,0),
     (-4,1,7,-14,7,1/2,7/2,7),(3/2,-3/2,3/2,7,-20,11,9/2,-7/2),(-17/2,11,-3/2,1/2,11,-16,19/2,-3/2),
     (7,-9/2,4,7/2,9/2,19/2,-14,1),(-9/2,0,0,7,-7/2,-3/2,1,-14)] 

\textcolor{magenta}{# type (2,8,5)}
u = [x^2*v^2,x^2*v*w,x^2*w^2,x*y*v^2,x*y*v*w,x*y*w^2,y^2*v^2,y^2*v*w,y^2*w^2]
R = [(0,1,0,0,1,0,1,0,1),(1,1,1,1,0,1,0,0,1),(1,1,0,1,1,1,1,0,1),(1,1,1,1,1,1,1,1,1),(1,0,1,0,1,0,1,0,0),
     (0,0,1,0,1,0,0,1,0),(0,0,0,0,0,0,1,1,1)] 
Q = [(-8,-8,10,8,-10,10,-2),(-8,-56,-26,100,-8,-40,-20),(10,-26,-50,74,17,-61,1),(8,100,74,-208,9,99,27),
     (-10,-8,17,9,-18,20,11),(10,-40,-61,99,20,-78,-9),(-2,-20,1,27,11,-9,-6)] 

\textcolor{magenta}{# type (2,8,4)}
u = [x^2*v^2,x^2*v*w,x^2*w^2,x*y*v^2,x*y*v*w,x*y*w^2,y^2*v^2,y^2*v*w,y^2*w^2]
R = [(0,1,0,1,1,1,0,1,1),(1,0,0,0,1,1,0,1,1),(1,0,0,0,1,0,0,0,1),(0,1,1,0,0,0,1,1,0),(1,1,1,0,0,0,1,1,0),
     (0,1,1,1,0,1,1,1,0)] 
Q = [(-8,5,2,9,-11,2),(5,-4,-1,-5,5,0),(2,-1,0,-9,11,-1),(9,-5,-9,-26,24,-1),(-11,5,11,24,-26,5),(2,0,-1,-1,5,-4)] 
\end{Verbatim}
\end{table}

\begin{example}
Let ({\tt u},{\tt R},{\tt Q}) be the entry in \TAB{cert}, where 
the type is~$(2,6,4)$.
The list {\tt u} encodes the birational map $\mu\c \S^1\times\S^1\dto \P^6$
that sends $(x:y;v:w)$ to
$(x^2 v^2 + y^2 v^2: x^2 v w: x^2 w^2: x y v w: x y w^2: y^2 v w: y^2 w^2)$.
The matrix~{\tt R} defines a linear projection $\rho\c \P^6\dto \P^5$
that sends $x$ to
$(x_0 + x_3 + x_4 + x_5: x_0 + x_3 + x_5 + x_6: x_1 + x_2 + x_3 + x_5: x_4: 
 x_1 + x_3 + x_4 + x_5: x_3 + x_4 + x_5)$.
The quadratic form associated to the symmetric matrix~{\tt Q}  
defines a hyperquadric~$Q\subset\P^5$ 
so that
$
Q:=\set{x\in\P^5}{
-6 x_0^2 + 6 x_0 x_1 - 7 x_0 x_2 + 19 x_1 x_2 + 11 x_0 x_3 + 19 x_1 x_3 - 2 x_2 x_3 
- 27 x_3^2 + 28 x_0 x_4 - 40 x_1 x_4 - 28 x_3 x_4 - 12 x_4^2 - 15 x_0 x_5 + 15 x_1 x_5 
+ 9 x_2 x_5 + 33 x_3 x_5 + 3 x_4 x_5 - 6 x_5^2=0}$.
Since $Q$ has signature~$(4+1,1)$, there exists 
projective transformation
$\nu\c\P^{5}\to\P^{5}$ \st $\nu(Q)=\S^4$.
The Zariski closed image of~$\nu\circ\rho\circ\mu$ 
is a celestial surface of type~$(2,6,4)$ (see forward \LEM{uRQ}).
\END
\end{example}

\begin{lemma}
\label{lem:verify}
Let $(\mu,\rho,Q)$ be encoded by the entry ({\tt u},{\tt R},{\tt Q}) in \TAB{cert}
that is marked with type~$(\lambda,d,n)$.
\begin{Mlist}
\item The matrix {\tt Q} is symmetric and has signature $(n+1,1)$.
\item The rank of the $(n+2)\times(d+1)$ matrix {\tt R} is $n+2$.
\item If $q$ is the quadratic form associated to the matrix {\tt Q} then $q\circ\rho\circ\mu=0$.
\item The components of $\rho\circ\mu$ are $n+1$ linearly independent forms of bidegree~$(2,2)$.
\item If $d=8$, then the linear series of $\rho\circ\mu$ is base point free.
\item If $(\lambda,d)=(3,6)$, then the linear series of $\rho\circ\mu$ 
has base points
$(1:-\Mi;1:\Mi)$ and $(1:\Mi;1:-\Mi)$.
\item If $(\lambda,d)=(2,6)$, then the linear series of $\rho\circ\mu$ 
has base points
$(1:\Mi;1:0)$ and $(1:-\Mi;1:0)$.
\end{Mlist}
\end{lemma}

\begin{proof}
Straightforward verification.
See \citep[Algorithm~1]{nls-basepoint} for computing base points of a linear series.
\end{proof}

\begin{lemma}
\label{lem:uRQ}
Let $(\mu,\rho,Q)$
be encoded by an entry ({\tt u},{\tt R},{\tt Q}) in \TAB{cert}
that is marked with type~$(\lambda,d,n)$.
Suppose that 
$Z_N\subset\P^d$
and
$Z\subset\P^{n+1}$ are the Zariski closed images of the maps~$\mu$ and $\rho\circ\mu$,
\resp.
\begin{Mlist}

\item[\bf a)]
The map $\mu\c \S^1\times\S^1\dto Z_N$ is birational, $\deg Z_N=d$
and $Z_N$ is covered by $\lambda$ pencils of conics.

\item[\bf b)]
If $\nu\c\P^{n+1}\to\P^{n+1}$ is a projective transformation \st $\nu(Q)=\S^n$, then
the Zariski closed image of the bidegree~$(2,2)$
birational map~$\nu\circ\rho\circ\mu$ is 
a celestial surface of type~$(\lambda,d,n)$.
There exists a projective transformation~$\nu$ \st $\nu(Q)=\S^n$.

\item[\bf c)]
If $X$ is a celestial surface of type~$(\lambda,d,n)$,
then 
there exists a linear projection $\tilde{\rho}\c Z_N\to X$
\st $\tilde{\rho}\circ\mu\c \S^1\times\S^1\dto X$ is birational
and of bidegree~$(2,2)$.

\end{Mlist}
\end{lemma}

\begin{proof}
a)
First suppose that $d=6$.
It follows from \LEM{verify} that 
the linear series of~$\mu$ consists of linear independent bidegree~$(2,2)$ forms
that vanish at the complex conjugate base points $p_1$ and $p_2$.
These base points lie either in general position or 
are contained in a fiber of~$\pi_1\c \S^1\times\S^1\to\S^1$ (recall \RMK{S1S1}).
We illustrated the configurations in \FIG{C1C2}, 
where the vertical and horizontal line segments
correspond to the fibers of $\pi_1$ and $\pi_2$, \resp.
\begin{figure}[!ht]
\centering
\def\rad{0.3}
\newcommand{\disc}[2] {\draw[draw=black, fill=red!50] (#1,#2) circle [radius=\rad];}
\newcommand{\grid}{
\draw[draw=gray ] (0,3) -- (3,3);
\draw[draw=gray ] (0,2) -- (3,2);
\draw[draw=gray ] (0,1) -- (3,1);
\draw[draw=gray ] (0,0) -- (3,0);
\draw[draw=gray ] (0,3) -- (0,0);
\draw[draw=gray ] (1,3) -- (1,0);
\draw[draw=gray ] (2,3) -- (2,0);
\draw[draw=gray ] (3,3) -- (3,0);
}
\begin{tabular}{c@{\hspace{3cm}}c}
\begin{tikzpicture}[scale=0.4]\grid;\disc{1}{2};\disc{2}{1};\end{tikzpicture}&
\begin{tikzpicture}[scale=0.4]\grid;\disc{1}{2};\disc{1}{1};\end{tikzpicture}
\\[-1mm]
{\bf C1} & {\bf C2}
\end{tabular}
\caption{}
\label{fig:C1C2}
\end{figure}
Let $\tau\c Y\to \S^1\times\S^1$ be the birational morphism that blows up $p_1$ and $p_2$
and recall from \LEM{BGE}a that $-\k=2\,\l_0+2\,\l_1-\p_1-\p_2$ 
is the anticanonical class of~$Y$.
Notice that $\mu\circ \tau\c Y\to Z_N$ is the map associated to~$-\k$
and thus $Z_N$ is the anticanonical model of~$Y$.
It follows from $-\k$ being nef and big and from Reider's lemma that $\mu\circ \tau$
and thus $\mu$ itself is birational.
Since $-\k$ is not orthogonal to the class of a $(-1)$-curve,
we find that $\mu\circ \tau$ is a smooth model for~$Z_N$.
It follows from \LEM{verify} and \LEM{BGE}b that 
$B(Z_N)=\{b_{12}\}$ or $B(Z_N)=\emptyset$
if $(\lambda,d)=(2,6)$ and $(\lambda,d)=(3,6)$, \resp.
Therefore, by \LEM{BGE}c, we have that $G(Z_N)$
is equal to either $\{g_0,g_1,g_{12}\}$ or $\{g_0,g_1\}$ 
and we verify that $\lambda=|G(Z_N)|$.
Suppose that $C\subset Z_N$ is a conic so that $[C]\in G(Z_N)$.
It follows from the Riemann-Roch 
and Kawamata vanishing theorems that $h^0([C])=2$
and thus $C$ is a member of a pencil of conics.
Hence, we conclude that $Z_N$ is 
covered by $\lambda$ pencils of conics as asserted.
If $d=8$, then the proof is similar, 
except in this case the linear series of~$\mu$ is base point free,
$B(Z_N)=\emptyset$ and $G(Z_N)=\{g_0,g_1\}$.

b)
It follows from \LEM{verify}
that $Z_N$ and $Z$ are not contained in a hyperplane section.
Moreover, we verified that
the base points of the linear series associated to
$\rho\circ\mu$ and $\mu$ are the same.
Therefore the linear series of~$\mu$ is a completion of the linear series of~$\rho\circ\mu$,
or equivalently,
$Z_N$ is a linear normalization of~$Z$.
We also know from \LEM{verify} that $Z\subset Q$ and that $Q$
has the same signature as~$\S^n$.
Hence, there exists a projective 
transformation~$\nu\c \P^{n+1}\to \P^{n+1}$
\st $\nu(Q)=\S^n$ and $\nu(Z)\subset\S^n$.
We now conclude from a) that $\nu(Z)$ is a celestial surface
of type $(\lambda,d,n)$.

c)
Let $X_N\subset\P^d$ be linear normalization of~$X$
and let $\varphi\c Y\to X_N$ be its smooth model.
Recall from \RMK{XN} that there exists 
a linear projection $\eta\c X_N\to X$
and notice that $\eta\circ\varphi\c Y\to X$ is a smooth 
model for~$X$.
First suppose that~$d=6$.
We know from \PRP{dp} that
$B(X)\subseteq\{b_{12}\}$, 
$G(X)\subseteq\{g_0,g_1,g_{12}\}$
and
there exists a birational morphism 
$\tau\c Y\to \S^1\times\S^1$ that contracts
two complex conjugate $(-1)$-curves 
to the complex conjugate blowup centers $p_1$ and $p_2$.
It follows from \LEM{BGE}[b,c] that  
$p_1$ and $p_2$ either lie in a general 
position or are contained in a fiber of~$\pi_1$.
Up to $\aut(\S^1\times\S^1)$ we may assume that one of the following two cases holds: 
\begin{Mlist}
\item 
$\lambda=3$,
$p_1=(1:-\Mi;1:\Mi)$ and $p_2=(1:\Mi;1:-\Mi)$ (see {\bf C1} in \FIG{C1C2}), or

\item 
$\lambda=2$,
$p_1=(1:\Mi;1:0)$ and $p_2=(1:-\Mi;1:0)$ (see {\bf C2} in \FIG{C1C2}).
\end{Mlist}
We know from \LEM{BGE}a that $-\k=2\,\l_0+2\,\l_1-\p_1-\p_2$ 
is the anticanonical class and thus the class of a hyperplane section of~$X_N$.
Therefore, the linear series of $\varphi\circ \tau^{-1}$
consist of all the bidegree~$(2,2)$ forms
that vanish at $p_1$ and $p_2$ with multiplicity one.
It follows from \LEM{verify}
that the components of~$\mu$ define a basis for such a linear series.
Therefore, there exists a projective isomorphism $\alpha\c Z_N\to X_N$
so that $\tilde{\rho}:=\eta\circ\alpha$
exists as is asserted.
If $d=8$, then the proof is similar, except in this
$Y\cong\S^1\times\S^1$ and $-\k=2\,\l_0+2\,\l_1$.
\end{proof}

\LEM{uRQ} suggests the following algorithm 
for computing examples of celestial surfaces. 
Each celestial surface
is reached by this algorithm,
except for Veronese surfaces and 2-circled Darboux cyclides
whose circle graphs are connected.

\begin{algorithm}[constructing examples of celestial surfaces]
\label{alg:construct}
~
\begin{itemize}[itemsep=0mm,topsep=-1mm, leftmargin=5mm]

\item {\bf Input:}
The set $B(X)$ and type $(\lambda,d,n)$
of some celestial surface~$X$ with smooth model~$Y\to X$
\st $Y$ is isomorphic to a blowup of $\S^1\times\S^1$.

\item {\bf Output:}
A birational map $\S^1\times\S^1\dto X$ of bidegree $(2,2)$
\st $X\subset\S^n$ is a celestial surface of type $(\lambda,d,n)$.

\item {\bf Method:} 
\begin{enumerate}[itemsep=0mm,topsep=0mm, leftmargin=0mm]

\item Use \LEM{BGE} to translate $B(X)$ into an explicit configuration of 
base points of a linear series (see \FIG{C1C2} and \LEM{D} for examples). 

\item
Apply 
\citep[Algorithm~2]{nls-basepoint} to construct a linear 
series of bidegree $(2,2)$ forms that pass through the base points with 
multiplicity one
and let $\mu\c \S^1\times\S^1\dto \P^d$ be 
the map associated to this linear series.

\item
Let $\rho\c \P^d\dto \P^{n+1}$ be a random linear projection.
Compute the graded ideal $I(Z)\subset\Q[y_0,\ldots,y_d]$ 
of the image $Z\subset\P^{n+1}$ of $\rho\circ\mu$
and let $(q_i)_i$ be the generators of the vector space of quadratic forms in $I(Z)$.
Let
$q:=\sum_ic_i\,q_i$,
where
$c_i\in\Z$ are random coefficients.
We compute $\deg Z$ and $\dim Z$ via the Hilbert polynomial of~$I(Z)$.

\item 
If either $(\deg Z,\dim Z)\neq (d,2)$ or if the signature 
of~$q$ is not equal to~$(n+1,1)$, then go back to and repeat step~3.
Otherwise, we diagonalize the symmetric matrix~{\tt Q} associated to~$q$
and obtain
a projective isomorphism $\nu\c\P^{n+1}\to\P^{n+1}$
\st $\nu(Q)=\S^n$, 
where $Q\subset\P^{n+1}$ is defined as the zero-set of~$q$.

\item 
We output the composition $\nu\circ\rho\circ\mu\c\S^1\times\S^1\dto X$. \END
\end{enumerate}
\end{itemize}
\end{algorithm}
 
\ALG{construct} has been implemented in \citep[{\tt orbital}]{orbital}
and was used for computing the examples in \TAB{cert}.
An open problem is compute $\rho$ and $q$ in step~3
\st each component of the projective transformation~$\nu$ has small coefficients in say $\{-1,0,1\}$.

A Veronese surface of type $(\infty,4,4)$
can be constructed with the same method as \ALG{construct},
but instead we set $\mu\c \P^2\to\P^5$ 
to be a map whose components generate the vector space
of quadratic forms on $\P^2$.


\begin{proposition}
\label{prp:type}
Each of the following types is realized by some celestial surface:
$(2,8,n)$ for $3\leq n\leq 7$, $(3,6,5)$, $(3,6,4)$, $(2,6,5)$ $(2,6,4)$
and $(\infty,4,4)$.
\end{proposition}

\begin{proof}
See \citep[(23.6), (23.7)]{kol2} for implicit equations
of a celestial surface of type $(\infty,4,4)$.
The remaining cases are now a direct consequence of \LEM{uRQ}b.
\end{proof}

\begin{remark}
Suppose that $\nu\circ\rho\circ\mu\c\S^1\times\S^1\dto X$ is the output of \ALG{construct}
with input $B(X)=\{b_{12}\}$ and $(2,6,4)$.
Let $q\in\S^1\times\S^1$ be a general point that lies in the same fiber as 
the complex conjugate base points (see {\bf C2} at \FIG{C1C2}).
It follows from \LEM{B}
that this fiber is contracted by $\nu\circ\rho\circ\mu$
so that $(\nu\circ\rho\circ\mu)(q)$ is a nodal singularity of $X\subset\S^4$. 
If we choose the center of stereographic projection~$\pi$ on this
singular point in~$X(\R)$, then $\pi(X(\R))\subset \R^4$ 
is covered by a pencil of lines and a pencil of circles.
A particularly nice example of such a surface was provided by 
Mikhail Skopenkov during private communication:
\[
\R^2\to\R^4,\quad (u,v)\mapsto (u\,v^2 - v,~ u\,v + v^2,~ u\,v + 1, u - v)/(v^2 + 1).
\]
If $d=8$, then
an alternative parametrization for the real points $X(\R)\subset S^7$ 
of a surface of type $(2,8,7)$ is as follows with $0\leq \alpha,\beta\leq 2\pi$:
\[
(\cos\alpha,\sin\alpha,\cos\beta,\sin\beta,
\cos(\alpha+\beta),\sin(\alpha+\beta),\cos(\alpha-\beta),\sin(\alpha-\beta))/2.
\]
The surface of type $(2,8,3)$ in \FIG{intro} was constructed as 
a stereographic projection
$\pi(\set{c \star d\in S^3}{a\in A,~b\in B})$, 
for some circles $A,B\subset S^3$, where~$\star$ denotes the Hamiltonian product for unit quaternions.
\END
\end{remark}


\section{Impossible types}
\label{sec:type}

In this section we show that a celestial surface cannot have type $(3,6,3)$, $(2,6,3)$ or $(\infty,4,3)$
so that a surface that contains $\lambda\geq 3$ circles through
a general point must be of type $(\infty,2,2)$, $(\infty,4,4)$, $(3,6,5)$, $(3,6,4)$ 
or $(\lambda,4,3)$ by \PRP{dp}. 
It will follow from \SEC{dar} that $\lambda\leq 6$
so that \THM{T} holds.

\begin{definition}
\label{def:P}
Suppose that $Z\subset \R^3$ is a surface. 
We denote by $\P(Z)\subset\P^3$ the Zariski closure of $\iota(Z)\subset \P^3$,
where $\iota\c \R^3\hookrightarrow \P^3$ sends $(x_1,x_2,x_3)$ 
to the projective point $(1:x_1:x_2:x_3)$.
\END
\end{definition}

The \df{delta invariant} $\delta_p(C)$ of a complex point $p$ in a
planar curve $C\subset\P^2$
is defined as $\Sigma_{q\in \cI_p}m_q(m_q-1)/2$,
where $\cI_p$ is the set that contains $p$ and complex points that are infinitely near to $p$,
and $m_q$ denotes the multiplicity of a strict transform of $C$ at~$q$ 
(see \DEF{blowup} and \citep[Example~V.3.9.2]{har}). 
Informally, we may think of~$\delta_p(C)$
as the number of double points that are concentrated at $p$ (see \citep[page 85]{delta}).
Notice that $\delta_p(C)>0$ \Iff $p\in \sng C$, and that 
$\delta_p(C)=\delta_{\sigma(p)}(C)$.

\begin{lemma}
\label{lem:s}
If $Z\subset\R^3$ is a $\lambda$-circled celestial surface of degree $d$
and if $H$ is a general hyperplane section of $\P(Z)$, then
$\sum_{p\in H}\delta_p(H)=\tfrac{1}{2}(d-1)(d-2)-s$,
where $s=0$ or $s=1$ if $\lambda=\infty$ and $\lambda<\infty$, \resp.
\end{lemma}

\begin{proof}
Since $H$ is a planar curve of degree~$d$
it follows from \citep[Section~2.4.6]{intersection} that
the geometric genus of $H$ is equal to
$p_g(H)=\frac{1}{2}(d-1)(d-1)-\sum_{p\in H}\delta_p(H)$.
It follows from \citep[Theorem~5 and Theorem~8]{conical}
that the sectional genus $p_g(H)$ is as asserted.
Alternatively, notice that
$p_g(H)=\tfrac{1}{2}([H]^2+[H]\cdot\k)+1$
by \citep[Proposition~IV.1.1 and Exercise~V.1.3]{har}
and it follows from \THM{mc}
that $[H]=-\frac{2}{3}\k$ 
or $[H]=-\k$
if $\lambda=\infty$ and $\lambda<\infty$, \resp.
\end{proof}

\begin{lemma}
\label{lem:delta}
If $X\subset\P^n$ is a surface
and $H\subset X$ a general hyperplane section,
then 
$p\in\sng H$ \Iff $p\in H\cap \sng X$.
\end{lemma}

\begin{proof}
First suppose that $p\in \sng H$
and let $\varphi\c Y\to X$ denote the smooth model of~$X$.
By Bertini's theorem \citep[page~137]{gh}
the preimage $\varphi^{-1}(H)$ for the general hyperplane section $H$ is smooth. 
Since $\varphi$ is an isomorphism outside the preimage of singular locus of~$X$,
it follows that $H$ is smooth outside singular locus of~$X$
and thus $p\in \sng X$ as asserted.
Now suppose that $p\in H\cap \sng X$.
Recall from \RMK{XN} that
there exists a linear projection $\eta\c X_N\to X$ from the linear normalization~$X_N$.
Let $H_N:=\eta^{-1}(H)$, $V_N:=\eta^{-1}(\sng X)$ and
$p_N:=\eta^{-1}(p)$ so that $p_N\subseteq V_N\cap H_N$.
Since $H$ is general, it follows that 
either $H_N$ is already singular at $p_N$ or $|p_N|>1$.
If $|p_N|>1$, then the local complex analytic branches of~$H_N$ with respective centers in~$p_N$ 
are mapped to branches of~$H$ each centered at~$p$.
Thus $H$ must be singular at~$p$ as asserted.
\end{proof}

\begin{lemma}
\label{lem:notype}
If $Z\subset\R^n$ is a celestial surface of type $(\lambda,d,n)$,
then $(\lambda,d,n)\notin\{(\infty,4,3),(2,6,3),(3,6,3)\}$.
\end{lemma}

\begin{proof}
Suppose by contradiction that $(d,n)=(6,3)$ with $\lambda<\infty$.
We may assume after a M\"obius transformation that $\deg Z=d$
so that $Z$ is compact in $\R^3$.
Thus there exists a general
hyperplane section $H$ of $\P(Z)$ \st $H\bas{\R}=\emptyset$.
We apply \LEM{s} and \LEM{delta} and find that $\sum_{p\in H\,\cap\,\sng\P(Z)}\delta_p(H)=9$.
If $p,q\in H$ are complex conjugate points, then $\delta_p(H)=\delta_q(H)$
as a direct consequence of the definitions.
We arrived at a contradiction as $\sum_{p\in H\,\cap\,\sng\P(Z)}\delta_p(H)$ 
must be even.
To show that $(\lambda,d,n)\neq(\infty,4,3)$
we use the same argument as before, but instead 
$\sum_{p\in H\,\cap\,\sng\P(Z)}\delta_p(H)=3$.
\end{proof}

We are now ready to answer the question in \citep[Section~5]{pot2}.

\begin{proof}[Proof of \COR{3}.]
Suppose that $Z\subset\R^3$ contains $\lambda\geq 3$ circles through a general point.
We know from \PRP{dp} (or alternatively from \LEM{even} and 
\citep[Theorems~5--8 and Proposition~1]{conical}) 
that $Z$ has type either 
$(\infty,2,2)$, $(\lambda,4,3)$, $(2,6,3)$, $(3,6,3)$ or $(\infty,4,3)$.
The proof is now concluded by \LEM{notype}
as the latter three types are impossible.
\end{proof}

\section{Darboux cyclides}
\label{sec:dar}

In this section we classify celestial Darboux cyclides
in order to complete the proof of \THM{B} and its corollaries.
The geometry of a celestial Darboux cyclide $X\subset\S^3$
is to a large extend determined by its N\'eron-Severi lattice~$N(X)$.
We classify $N(X)$
by classifying unimodular involutions $\sigma_*$ 
and subsets $B(X)$
up to $\aut N(X)$.

\begin{remark}
\COR{graph} can be approached 
using the methods from \cite{cool,pot2,tak2}
and it would be interesting to see such an alternative proof.
\END
\end{remark}

\begin{lemma}
\label{lem:sigma}
If $X\subset\S^3$ is a celestial Darboux cyclide with smooth model~$Y\to X$, 
then up to $\aut N(X)$ the unimodular involution $\sigma_*$ and $Y$
are characterized by a row of \TAB{sigma},
where the centers of blow up are in sufficiently general position.
Moreover, if $\sigma_*\in\{2A_1,3A_1,D_4\}$, then 
$\set{g\in G(X)}{\sigma_*(g)=g}$ is contained in 
$\{g_0,g_1,g_2,g_3,g_{12},g_{34}\}$, $\{g_{12},g_{34}\}$ and $\{g_1,g_2\}$, \resp.
If $X$ is in addition smooth and $\sigma_*$ is $3A_1$ or $D_4$, 
then $X$ is a $S1$ cyclide and $S2$ cyclide, \resp. 
\end{lemma}

\begin{table}[!ht]
\caption{See \LEM{sigma}.}
\label{tab:sigma}
\centering
\begin{tabular}{l@{\hspace{1cm}}l}
$\sigma_*$ & smooth model                                               \\\hline
$2A_1$     & $\S^1\times\S^1$ blown up in two pairs of complex conjugate points \\
$3A_1$     & $\S^2$ blown up in two pairs of complex conjugate points           \\
$D_4$      & $X$                                                                    
\end{tabular}
\end{table}

\begin{proof}
Recall from \RMK{dp} and \PRP{dp} that $Y$ is a weak del Pezzo surface.
As $\sigma_*\c N(X)\to N(X)$ 
is induced by the real structure it leaves the canonical class~$\k$ of~$X$ invariant.
The matrix defining $\sigma_*$ is an involution and thus has eigenvalues~$\pm 1$.
Notice that $\sigma_*(v)\cdot \sigma_*(\k)=-v\cdot k$ for all $v$
in the eigenspace of $-1$. Hence this eigenspace
is contained in the following inner product space 
$V_{\k}(X):=(\set{ c\in N(X) }{ \k\cdot c=0 }\otimes_\Z\R,~ \cdot)$.
By \citep[8.6.3]{dol} the set $R(X):=\set{c\in V_{\k}(X)}{c^2=-2}$ forms a root system
of Dynkin type~$D_5$ in $V_{\k}(X)$.
The intersection of a root system with a subspace is a root subsystem
and thus $S_\sigma(X):=\set{ c\in R(X) }{ \sigma_*(c)=-c }$ forms a root subsystem of $R(X)$.
We know from \citep[Corollary~2.1]{wal1} that $S_\sigma(X)$ has Dynkin type either 
\[
A_0:=\emptyset,\qquad  
A_1,\qquad               
2A_1,\qquad              
2A_1',\qquad           
3A_1
\qquad\text{or}\qquad
D_4.
\]
Root subsystems with the same Dynkin type $2A_1$ may not
be the same up to $\aut N(X)$ and therefore one is marked with $'$.

A linearly independent subset $\{s_1,\ldots,s_r\}\subset S_\sigma(X)$ 
is called a \df{root base} if $s_i\cdot s_j\geq 0$ for all $1\leq i,j\leq r$.  
The Dynkin type of $S_\sigma(X)$ is defined by the type of the incidence diagram
of a root base and does not depend on the choice.
We may assume \Wlog that $s_i\sim\{b_0,b_1,b_{12}\}$ (see \citep[8.2.3]{dol}).

Let 
$G'(X):=\set{c\in N(X)}{\sigma_*(c)=c,~c\sim \{g_0,g_2,g_{12}\}}$ and
$E'(X):=\set{c\in N(X)}{\sigma_*(c)=c,~c\sim \{e_1,e_1',e_{01}\}}$.
If $E'(X)>0$, then it is a straightforward consequence of \citep[Lemma~8.2.22]{dol} that
$\set{c\in E(X)}{\sigma_*(c)=c}>0$.
Thus by \LEM{BGE} and \LEM{even}, we require  
that $|G'(X)|\geq 2$ and $|E'(X)|=0$.

For each of the six Dynkin types of $S_\sigma(X)$ 
we find a basis of given type and 
from this basis we construct 
explicit coordinates for~$\sigma_*$. 
We require that $|G'(X)|\geq 2$ and $|E'(X)|=0$ for the constructed $\sigma_*$.

For example, suppose that $S_\sigma(X)$ has Dynkin type $D_4$.
Up to $\aut N(X)$, $S_\sigma(X)$ has root base
$\{\l_1-\p_1-\p_2$, $\p_1-\p_2$, $\p_2-\p_3$, $\p_3-\p_4\}$.
We consider the following three matrices:
\[
\setlength{\abovedisplayskip}{10pt}\setlength{\belowdisplayskip}{10pt}
B=
\left(\begin{smallmatrix}
 0 & 0 & 0 & 0 \\
 1 & 0 & 0 & 0 \\
-1 & 1 & 0 & 0 \\
-1 &-1 & 1 & 0 \\
 0 & 0 &-1 & 1 \\
 0 & 0 & 0 &-1 \\
\end{smallmatrix}\right)
,~
J=
\left(\begin{smallmatrix}
 0 & 1 & 0 & 0 & 0 & 0\\
 1 & 0 & 0 & 0 & 0 & 0\\
 0 & 0 &-1 & 0 & 0 & 0\\
 0 & 0 & 0 &-1 & 0 & 0\\
 0 & 0 & 0 & 0 &-1 & 0\\
 0 & 0 & 0 & 0 & 0 &-1\\
\end{smallmatrix}\right)
,~
V=
\left(\begin{smallmatrix}
 0 & 0 & 0 & 0 & 2 & 0\\
 1 & 0 & 0 & 0 & 0 & 1\\
-1 & 1 & 0 & 0 &-1 & 0\\
-1 &-1 & 1 & 0 &-1 & 0\\
 0 & 0 &-1 & 1 &-1 & 0\\
 0 & 0 & 0 &-1 &-1 & 0\\
\end{smallmatrix}\right)
.
\]
The four colums of $B$ correspond to the generators of this root base
and are eigenvectors for the eigenvalue $-1$ of $\sigma_*$.
The matrix $J$ defines the intersection product for $N(X)$.
The matrix $V$ is obtained by augmenting $B$ with two column vectors that generate the kernel of the matrix $B^\top\cdot {J}$.
The appended two columns are eigenvectors for the eigenvalue $1$ of $\sigma_*$.
Let $D$ be the diagonal matrix with eigenvalues $(-1,-1,-1,-1,1,1)$ on the diagonal.
We denote the matrix corresponding to $\sigma_*$ by $M$.
Since $M\cdot V=V\cdot D$, it follows that $V\cdot D\cdot V^{-1}$ is the matrix corresponding to 
involution $\sigma_*$ \st $\sigma_*(\l_0)=g_3$, 
$\sigma_*(\l_1)=\l_1$ and $\sigma_*(\p_i)=\l_1-\p_i$ for $1\leq i\leq 4$.
We verify that $|E'(X)|=0$ and find that $G'(X)=\{g_1,g_2\}$ so that $G'(X)=G(X)$. 
Notice that $\sigma_*$ is unique up to inner automorphisms in $\aut N(X)$
and that we used the same notation~$D_4$ in \SEC{intro}.
This concludes the proof of this lemma for the case that
$S_\sigma(X)$ has Dynkin type $D_4$. 

The remaining cases are similar and we verify that $|E'(X)|=0$ and $|G'(X)|\geq 2$
only if $S_\sigma(X)$ has Dynkin type either $2A_1$, $3A_1$ or $D_4$.
The characterizations of~$Y$ is now a direct consequence of \RMK{S1S1} 
and \PRP{dp}.
The topological characterizations as S1 cyclide or S2 cyclide follows from \citep[Corollary~3.2]{wal1}.
\end{proof}

\begin{lemma}
\label{lem:A3D4}
If $X\subset\S^3$ is a celestial Darboux cyclide 
and either $B(X)=\emptyset$ or $\sigma_*\neq 2A_1$, 
then $\sigma_*$, $B(X)$, $\sng X$ and $G(X)$
are up to $\aut N(X)$ characterized by a row in \TAB{A3D4}. 
We overlined a class $g$ in $G(X)$ if $\sigma_*(g)=g$.
\end{lemma}

\begin{table}[!ht]
\caption{See \LEM{A3D4}.}
\label{tab:A3D4}
\centering
\begin{tabular}{l@{\hspace{1cm}}l@{\hspace{1cm}}l@{\hspace{1cm}}l}
$\sigma_*$ & $B(X)$               & $\sng X$ & $G(X)$ \\\hline
$2A_1$     & $\emptyset$          & $\AN$    & $\{\g_0,\g_1,\g_2,\g_3,\g_{12},\g_{34},g_{13},g_{24},g_{14},g_{23}\}$\\
$3A_1$     & $\{b_0\}$            & $\uAi$   & $\{g_0,g_1,\g_{12},\g_{34},g_{13},g_{24},g_{14},g_{23}\}$ \\
$3A_1$     & $\{b_{13},b_{24}'\}$ & $\uAii$  & $\{g_0,g_1,\g_{12},\g_{34},g_{14},g_{23}\}$ \\
$3A_1$     & $\emptyset$          & $\AN$    & $\{g_0,g_1,g_2,g_3,\g_{12},\g_{34},g_{13},g_{24},g_{14},g_{23}\}$ \\
$D_4$      & $\emptyset$          & $\AN$    & $\{
\g_1,
\g_2,
g_0,
g_3,
g_{12},
g_{34},
g_{13},
g_{24},
g_{14},
g_{23}
\}$ \\           
\end{tabular}
\end{table}

\begin{proof}
By \LEM{sigma} we have $\sigma_*\in\{2A_1,3A_1,D_4\}$ 
and we verify for each case the assertion using \LEM{BGE}. 
For example, if $\sigma_*$ is $D_4$, then $G(X)=\{g_1,g_2\}$ by \LEM{sigma},
$\sigma_*(\p_i-\p_j)=\p_j-\p_i$,
$\sigma_*(b_{ij})=-b_{ij}'$
and $g_2\cdot b_{ij}, g_2\cdot b_0<0$ for all $1\leq i < j\leq 4$
so that $B(X)=\emptyset$ by \LEM{BGE}.
The remaining cases are similar and thus straightforward.
\end{proof}

\begin{remark}
\label{rmk:D}
Suppose that $X\subset\S^3$ is a celestial Darboux cyclide with smooth model~$Y\to X$
\st $\sigma_*$ is $2A_1$.
Thus $Y$ is by \LEM{sigma} isomorphic to~$\S^1\times\S^1$ blown up in two pairs of complex conjugate points
$p_1$, $p_2$, $p_3$ and~$p_4$.
We visualize $\S^1\times\S^1$ as a square such that
horizontal and vertical line segments in the square correspond to the fibers of $\pi_1$ and $\pi_2$, \resp.
The real structure acts on this square by sending horizontal line segments to 
horizontal line segments and vertical line segments to vertical line segments.
A pair of complex conjugate blowup centers are depicted by either two squares or two circles.
We consider in \LEM{D} below the following configurations, where 
the four points do not lie on a curve of bidegree (1,1).
\begin{center}
\def\rad{0.3}
\newcommand{\disc}[2] {\draw[draw=black, fill=red!50] (#1,#2) circle [radius=\rad];}
\newcommand{\boxx}[2] {\draw[draw=black, fill=blue!50] (#1-\rad,#2-\rad) rectangle (#1+\rad,#2+\rad);}
\newcommand{\grid}{
\draw[draw=gray ] (0,3) -- (3,3);
\draw[draw=gray ] (0,2) -- (3,2);
\draw[draw=gray ] (0,1) -- (3,1);
\draw[draw=gray ] (0,0) -- (3,0);
\draw[draw=gray ] (0,3) -- (0,0);
\draw[draw=gray ] (1,3) -- (1,0);
\draw[draw=gray ] (2,3) -- (2,0);
\draw[draw=gray ] (3,3) -- (3,0);
}
\begin{tabular}{@{}c@{\hspace{5mm}}c@{\hspace{5mm}}c@{\hspace{5mm}}c@{\hspace{5mm}}c@{\hspace{5mm}}c@{\hspace{5mm}}c@{}}
\begin{tikzpicture}[scale=0.4]\grid;\disc{0}{3};\disc{1}{3};\boxx{2}{1};\boxx{3}{0};\end{tikzpicture}&
\begin{tikzpicture}[scale=0.4]\grid;\disc{0}{3};\disc{1}{3};\boxx{2}{0};\boxx{3}{0};\end{tikzpicture}&
\begin{tikzpicture}[scale=0.4]\grid;\disc{0}{3};\boxx{1}{3};\disc{2}{0};\boxx{3}{0};\end{tikzpicture}&
\begin{tikzpicture}[scale=0.4]\grid;\disc{0}{3};\disc{0}{0};\boxx{3}{3};\boxx{2}{0};\end{tikzpicture}&
\begin{tikzpicture}[scale=0.4]\grid;\disc{0}{3};\disc{1}{3};\boxx{3}{1};\boxx{3}{0};\end{tikzpicture}&
\begin{tikzpicture}[scale=0.4]\grid;\disc{0}{0};\disc{3}{0};\boxx{0}{3};\boxx{3}{3};\end{tikzpicture}&
\begin{tikzpicture}[scale=0.4]\grid;\disc{0}{0};\boxx{0}{3};\disc{3}{3};\boxx{3}{0};\end{tikzpicture}
\\[-3mm]{\bf D1} & {\bf D2  } & {\bf D3} & {\bf D4} & {\bf D5} & {\bf D6} & {\bf D7}\\
\end{tabular}
\end{center}
We denote an infinitely near blowup center as an overlapping square. 
If $p_2$ lies additionally on the pullback of a fiber containing $p_1$, 
then the overlapping square lies in the horizontal or vertical direction of the blowup center,  
like in {\bf E3} and {\bf E4}.
\begin{center}
\def\rr{0.2}
\def\rs{0.3}
\newcommand{\infN}[2] {
\draw[draw=black, fill=red!50] (#1,#2) circle [radius=\rs];
\draw[draw=black, fill=blue!50] (#1-\rr,#2-\rr+\rs) rectangle (#1+\rr,#2+\rr+\rs);
}
\newcommand{\infW}[2] {
\draw[draw=black, fill=red!50] (#1,#2) circle [radius=\rs];
\draw[draw=black, fill=blue!50] (#1-\rr-\rs,#2-\rr) rectangle (#1+\rr-\rs,#2+\rr);
}
\newcommand{\infS}[2] {
\draw[draw=black, fill=red!50] (#1,#2) circle [radius=\rs];
\draw[draw=black, fill=blue!50] (#1-\rr,#2-\rr-\rs) rectangle (#1+\rr,#2+\rr-\rs);
}
\newcommand{\infSE}[2] {
\draw[draw=black, fill=red!50] (#1,#2) circle [radius=\rs];
\draw[draw=black, fill=blue!50] (#1-\rr+\rs,#2-\rr-\rs) rectangle (#1+\rr+\rs,#2+\rr-\rs);
}
\newcommand{\infNW}[2] {
\draw[draw=black, fill=red!50] (#1,#2) circle [radius=\rs];
\draw[draw=black, fill=blue!50] (#1-\rr-\rs,#2-\rr+\rs) rectangle (#1+\rr-\rs,#2+\rr+\rs);
}
\newcommand{\infSW}[2] {
\draw[draw=black, fill=red!50] (#1,#2) circle [radius=\rs];
\draw[draw=black, fill=blue!50] (#1-\rr-\rs,#2-\rr-\rs) rectangle (#1+\rr-\rs,#2+\rr-\rs);
}
\newcommand{\grid}{
\draw[draw=white, fill=white] (0,0) circle [radius=\rs];
\draw[draw=white, fill=white] (0,3) circle [radius=\rs];
\draw[draw=gray ] (0,3) -- (3,3);
\draw[draw=gray ] (0,2) -- (3,2);
\draw[draw=gray ] (0,1) -- (3,1);
\draw[draw=gray ] (0,0) -- (3,0);
\draw[draw=gray ] (0,3) -- (0,0);
\draw[draw=gray ] (1,3) -- (1,0);
\draw[draw=gray ] (2,3) -- (2,0);
\draw[draw=gray ] (3,3) -- (3,0);
}
\begin{tabular}{@{}c@{\hspace{10mm}}c@{\hspace{10mm}}c@{\hspace{10mm}}c@{}}
\begin{tikzpicture}[scale=0.4]\grid;\infSE{0}{3};\infNW{3}{0}\end{tikzpicture}&
\begin{tikzpicture}[scale=0.4]\grid;\grid;\infSE{0}{3};\infSW{3}{3}\end{tikzpicture}&
\begin{tikzpicture}[scale=0.4]\grid;\infS{0}{3};\infS{3}{3}\end{tikzpicture}&
\begin{tikzpicture}[scale=0.4]\grid;\infS{0}{3};\infN{3}{0}\end{tikzpicture}
\\[-3mm]{\bf E1} & {\bf E2} & {\bf E3} & {\bf E4} \\
\end{tabular}
\end{center}
It follows from \LEM{BGE} and \LEM{B} that 
a fiber in $Y\cong\S^1\times\S^1$ that contains two centers of blowup
is contracted via the smooth model~$Y\to X$ to a complex isolated singularity. 
\END
\end{remark}

\begin{lemma}
\label{lem:D}
If $X\subset\S^3$ is a celestial Darboux cyclide \st $B(X)\neq\emptyset$ 
and $\sigma_*=2A_1$, then $B(X)$, $\sng X$ and $G(X)$
are up to $\aut N(X)$ characterized by a row in \TAB{D}. 
Moreover, we
included configurations for the centers of blowup from \RMK{D}
and overlined classes in $G(X)$ \st $\sigma_*(g)=g$.
\end{lemma}

\begin{table}[!ht]
\caption{See \LEM{D}.}
\label{tab:D}
\centering
\begin{tabular}{llll}
$B(X)$                               & $\sng X$      & $G(X)$                                           &\\
\hline
$\{b_1,b_2\}$                        & $2\Ai$        & $\{\g_0,\g_1,\g_2,\g_3,\g_{12},g_{13},g_{24}\}$  &{\bf D3/E1}\\
$\{b_{13},b_{24},b_{14}',b_{23}'\}$  & $4\Ai$        & $\{\g_0,\g_1,\g_2,\g_3\}$                        &{\bf D7/E4}\\
$\{b_{12}\}$                         & $\uAi$        & $\{\g_0,\g_1,\g_2,g_{13},g_{24},g_{14},g_{23}\}$ &{\bf D1}\\
$\{b_{13},b_{24},b_{12}'\}$          & $\uAi+2\Ai$   & $\{\g_0,\g_1,\g_{34},g_{14},g_{23}\}$            &{\bf D4}\\  
$\{b_{12},b_{34}'\}$                 & $\uAii$       & $\{\g_0,\g_1,g_{13},g_{24},g_{14},g_{23}\}$      &{\bf D5}\\
$\{b_1,b_2,b_{12}\}$                 & $\uAiii$      & $\{\g_0,\g_1,\g_3,g_{13},g_{24}\}$               &{\bf E2}\\
$\{b_1,b_2,b_{12},b_{13}',b_{24}'\}$ & $\uAiii+2\Ai$ & $\{\g_0,\g_1\}$                                  &{\bf E3}\\
$\{b_{12},b_{34}\}$                  & $2\uAi$       & $\{\g_0,\g_1,\g_3,g_{13},g_{24},g_{14},g_{23}\}$ &{\bf D2}\\
$\{b_{12},b_{34},b_{13}',b_{24}'\}$  & $2\uAi+2\Ai$  & $\{\g_0,\g_1,g_{14},g_{23}\}$                    &{\bf D6}\\
\end{tabular}           
\end{table}

\begin{proof}
We follow the notation of \RMK{D}.
Since $B(X)\neq \emptyset$ it follows from \LEM{BGE}
that the four centers of blowup do not lie in general position,
but only in sufficiently general position.
We first consider the case that 
the four base points do not lie on a curve of bidegree~$(1,1)$.

For each value for 
$\alpha:=|\pi_1(P)|$
and 
$\beta:=|\pi_2(P)|$
with $P:=\{p_1,p_2,p_3,p_4\}$,
we list all possible
configurations up to symmetry by using the following three restrictions:
base points must come in complex conjugate pairs,
fibers that contain two non-conjugate base points come in pairs, and
at most two non-infinitely near base points are allowed to lie in the same fiber.
Observe that $2\leq \alpha,\beta\leq 4$, $(\alpha,\beta)\neq (4,4)$
and $(\alpha,\beta)=(\beta,\alpha)$ up to symmetry.
Thus we obtain the following table for non-infinitely near base points:
\begin{center}
\begin{tabular}{r|c@{\hspace{1cm}}c@{\hspace{1cm}}c@{\hspace{1cm}}c@{\hspace{1cm}}c}
$(\alpha,\beta)$ & (3,4)    & (2,4)          & (3,3) & (2,3) & (2,2) \\
configuration    & {\bf D1} & {\bf D2,D3} & {\bf D5} & {\bf D4}     & {\bf D6,D7}\\
\end{tabular}
\end{center}
For infinitely near base points we find that {\bf E1--E4} are all possible configurations,
as infinity near base points must come in complex conjugate pairs.
It follows from \PRP{C} that we need to verify
using \LEM{BGE} that 
$|\set{g\in G(X)}{\sigma_*(g)=g}|\geq 2$ for each configuration.
For example, for  configuration {\bf D4} we have
$B(X)=\{b_{13},b_{24},b_{12}'\}$, 
$G(X)=\{\g_0,\g_1,\g_{34},g_{14},g_{23}\}$
and thus $\sng X$ is $\uAi+2A_1$ by \LEM{B}.
The other cases are similar.

Notice that both {\bf D7} and {\bf E4} 
are configurations \st $\sng X$ is $4A_1$.
These two cases are equivalent, in the sense
that there exists an isomorphism between their associated N\'eron-Severi lattices.
For the first configuration we have that 
$B(X)=\{b_{12},b_{24},b_{14}',b_{23}'\}$.
We now apply the isomorphism $\mu\c N(X)\cong N(X)$
that sends 
$\l_0$, $\l_1$, $\p_1$, $\p_2$, $\p_3$, $\p_4$
to
$g_3$,
$g_{34}$,
$\l_0-\p_4$,
$\l_0-\p_3$,
$\l_1-\p_2$,
$\l_1-\p_1$, \resp.
The image of $B(X)$ via this isomorphism corresponds to the configuration of {\bf E4}
as $\mu(B(X))=\{b_{13}',b_{24}',b_1,b_2\}$.
The configurations {\bf D3} and {\bf E1} are equivalent for similar reasons.

Finally, we include the configuration where
all four base points lie on a curve of bidegree (1,1).
We go through each of the configurations {\bf D1--D7} and {\bf E1--E4},
but now include $b_0$ as an additional element for $B(X)$. 
We find that $|G(X)|<2$, unless $B(X)=\{b_0\}$.
However, this case is the same as $B(X)=\{b_{12}\}$ up to $\aut N(X)$
and thus we concluded the proof of this lemma.
\end{proof}

\begin{definition}
\label{def:U}
The \df{hyperplane at infinity} $H_\infty\subset \P^3$
is the hyperplane that contains 
the \df{absolute conic}
$U:=\set{x\in\P^3}{x_0=x_1^2+x_2^2+x_3^2=0}$.
A \df{complex circle} in $\P^3$ is defined as an irreducible complex conic
that intersects $U$ either tangentially or in two complex points.
\END
\end{definition}

We included a proof of the following classically known proposition due to the lack of suitable reference.

\begin{proposition}
\label{prp:Q} 
If $Z\subset\R^3$ is irreducible quadric surface, 
then, up to Euclidean similarity, the following 
data is defined by exactly one row of \TAB{quad}:
\begin{Mlist}
\item the quadratic form $q(x)$ \st $\P(Z)=\set{x\in \P^3}{ q(x)=0 }$, 
\item the number of circles through a general point that are not lines,
\item the number of lines through a general point, and
\item $|G(X)|$, where $X\subset \S^3$ is the M\"obius model of $Z$.
\end{Mlist}
\end{proposition}

\begin{table}[!hp]
\caption{See \PRP{Q} and let
$\alpha,\beta,\gamma \in\R_{>0}$ with $\alpha\neq\beta$ and $\gamma\neq 1$.}
\label{tab:quad}
{
\begin{tabular}{lccccc}
$q(x)                                             $ & singular & \#circles  & \#lines & $|G(X)|$       & name             \\
\hline
$\alpha\, x_1^2 + \beta\, x_2^2 + x_3^2 - x_0^2   $ & no       & $2$        & $0$     & $8           $ & EE            \\
$\gamma\, x_1^2 + \gamma\, x_2^2 + x_3^2 - x_0^2  $ & no       & $1$        & $0$     & $5           $ & CE            \\
$x_1^2 + x_2^2 + x_3^2 - x_0^2                    $ & no       & $\infty$   & $0$     & $\infty      $ & $\S^2$        \\ 
$\alpha\, x_1^2 + \beta\, x_2^2 + x_3^2 + x_0^2   $ & no       & $0$        & $0$     & $8           $ & empty-1       \\
$\gamma\, x_1^2 + \gamma\, x_2^2 + x_3^2 + x_0^2  $ & no       & $0$        & $0$     & $5           $ & empty-2       \\ 
$x_1^2 + x_2^2 + x_3^2 + x_0^2                    $ & no       & $0$        & $0$     & $\infty      $ & empty-3       \\
$\alpha\, x_1^2 + \beta\, x_2^2 + x_3^2           $ & yes      & $0$        & $0$     & $7           $ & point-1       \\
$\gamma\, x_1^2 + \gamma\, x_2^2 + x_3^2          $ & yes      & $0$        & $0$     & $4           $ & point-2       \\
$x_1^2 + x_2^2 + x_3^2                            $ & yes      & $0$        & $0$     & $\infty      $ & point-3       \\
\hdashline
$\alpha\, x_1^2 + \beta\, x_2^2 - x_3^2 + x_0^2   $ & no       & $2$        & $0$     & $8           $ & EH2           \\
$\alpha\, x_1^2 + \alpha\, x_2^2 - x_3^2 + x_0^2  $ & no       & $1$        & $0$     & $5           $ & CH2           \\
$\alpha\, x_1^2 + \beta\, x_2^2 - x_3^2 - x_0^2   $ & no       & $2$        & $2$     & $8           $ & EH1           \\
$\alpha\, x_1^2 + \alpha\, x_2^2 - x_3^2 - x_0^2  $ & no       & $1$        & $2$     & $5           $ & CH1           \\
$\alpha\, x_1^2 + \beta\, x_2^2 - x_3^2           $ & yes      & $2$        & $1$     & $7           $ & EO            \\
$\alpha\, x_1^2 + \alpha\, x_2^2 - x_3^2          $ & yes      & $1$        & $1$     & $4           $ & CO            \\
\hdashline
$\gamma\, x_1^2 + x_2^2  + x_3x_0                 $ & no       & $2$        & $0$     & $6           $ & EP            \\
$\alpha\, x_1^2 + \alpha\, x_2^2 + x_3x_0         $ & no       & $1$        & $0$     & $3           $ & CP            \\
$\gamma\, x_1^2 + x_2^2 - x_0^2                   $ & yes      & $2$        & $1$     & $5           $ & EY            \\
$\alpha\, x_1^2 + \alpha\, x_2^2 - x_0^2          $ & yes      & $1$        & $1$     & $2           $ & CY            \\
$\gamma\, x_1^2 + x_2^2 + x_0^2                   $ & yes      & $0$        & $0$     & $5           $ & point-4       \\
$x_1^2 + x_2^2 + x_0^2                            $ & yes      & $0$        & $0$     & $2           $ & point-5       \\
\hdashline
$\alpha\, x_1^2 - x_2^2 + x_3x_0                  $ & no       & $0$        & $2$     & $6           $ & HP            \\
$\alpha\, x_1^2 - x_2^2 - x_0^2                   $ & yes      & $0$        & $1$     & $5           $ & HY            \\
\hdashline
$\alpha\, x_1^2 + x_0x_2                          $ & yes      & $0$        & $1$     & $1           $ & PY            \\
\end{tabular}
}%
\end{table}

\begin{proof}
Let $Q$ denote the projective quadric $\P(Z)$. 
The Euclidean similarities act via the embedding $\iota\c\R^3\hookrightarrow\P^3$ 
on the hyperplane at infinity $H_\infty$ as the orthogonal group while leaving 
the absolute conic $U$ invariant. 
It follows that we can diagonalize the quadratic form associated to $Q\cap H_\infty$
and thus we may assume up to rotations and translations that $q(x)$ is of the form
\begin{equation}
\label{eqn:q}
q(x)=a_1\, x_1^2+ a_2\, x_2^2 + a_3\, x_3^2 + x_0\,(b_0\,x_0+b_1\,x_1+b_2\,x_2+b_3\,x_3),
\end{equation}
for $a_1,a_2,a_3,b_1,b_2,b_3\in \R$.
We make a case distinction on the curve $Q\cap H_\infty$ 
and $U\cap Q$ which are represented below by a solid red curve and 
a dashed blue circle, \resp.
\begin{center}
\begin{tabular}{@{}cccccc@{}}
\begin{tikzpicture}[scale=0.6]
\def\xx{1.4}
\draw[white,thick] (0,1.4) to (1.8,-0.4);\draw[white,thick] (0,-1.4) to (1.8,0.4);
\draw[white,thick] (-\xx,\xx) to (\xx,-\xx);\draw[white,thick] (\xx,\xx) to (-\xx,-\xx);
\draw[draw=blue, very thick, densely dashed] (0,0) circle [radius=1cm];
\draw[red,thick] (0,0) circle [radius=1.02cm];
\end{tikzpicture}
&
\begin{tikzpicture}[scale=0.6]
\def\xx{1.4}
\draw[white,thick] (0,1.4) to (1.8,-0.4);\draw[white,thick] (0,-1.4) to (1.8,0.4);
\draw[white,thick] (-\xx,\xx) to (\xx,-\xx);\draw[white,thick] (\xx,\xx) to (-\xx,-\xx);
\draw[draw=blue, very thick, densely dashed] (0,0) circle [radius=1cm];
\draw[red,thick] (-\xx,0) to [out=90, in=90] (\xx,0) to [out=270, in=270] (-\xx,0);
\draw[black,fill=green] (-0.7, 0.7) circle [radius=1.5mm];
\draw[black,fill=green] (-0.7,-0.7) circle [radius=1.5mm];
\draw[black,fill=green] ( 0.7, 0.7) circle [radius=1.5mm];
\draw[black,fill=green] ( 0.7,-0.7) circle [radius=1.5mm];
\end{tikzpicture}
&
\begin{tikzpicture}[scale=0.6]
\def\xx{1.4}
\draw[white,thick] (0,1.4) to (1.8,-0.4);\draw[white,thick] (0,-1.4) to (1.8,0.4);
\draw[white,thick] (-\xx,\xx) to (\xx,-\xx);\draw[white,thick] (\xx,\xx) to (-\xx,-\xx);
\draw[draw=blue, very thick, densely dashed] (0,0) circle [radius=1cm];
\draw[red,thick] (-1,0) to [out=90, in=90] (1,0) to [out=270, in=270] (-1,0);
\draw[black,fill=green] (-1,0) circle [radius=1.5mm];
\draw[black,fill=green] ( 1,0) circle [radius=1.5mm];
\end{tikzpicture}
&
\begin{tikzpicture}[scale=0.6]
\def\xx{1.4}
\draw[white,thick] (0,1.4) to (1.8,-0.4);\draw[white,thick] (0,-1.4) to (1.8,0.4);
\draw[white,thick] (-\xx,\xx) to (\xx,-\xx);\draw[white,thick] (\xx,\xx) to (-\xx,-\xx);
\draw[draw=blue, very thick, densely dashed] (0,0) circle [radius=1cm];
\draw[red,thick] (-\xx,\xx) to (\xx,-\xx);\draw[red,thick] (\xx,\xx) to (-\xx,-\xx);
\draw[black,fill=green] (-0.7, 0.7) circle [radius=1.5mm];
\draw[black,fill=green] (-0.7,-0.7) circle [radius=1.5mm];
\draw[black,fill=green] ( 0.7, 0.7) circle [radius=1.5mm];
\draw[black,fill=green] ( 0.7,-0.7) circle [radius=1.5mm];
\end{tikzpicture}
&
\begin{tikzpicture}[scale=0.6]
\def\xx{1.4}
\draw[white,thick] (0,1.4) to (1.8,-0.4);\draw[white,thick] (0,-1.4) to (1.8,0.4);
\draw[white,thick] (-\xx,\xx) to (\xx,-\xx);\draw[white,thick] (\xx,\xx) to (-\xx,-\xx);
\draw[draw=blue, very thick, densely dashed] (0,0) circle [radius=1cm];
\draw[red,thick] (0,1.4) to (1.8,-0.4);\draw[red,thick] (0,-1.4) to (1.8,0.4);
\draw[black,fill=green] (0.7,0.7) circle [radius=1.5mm];
\draw[black,fill=green] (0.7,-0.7) circle [radius=1.5mm];
\end{tikzpicture}
&
\begin{tikzpicture}[scale=0.6]
\def\xx{1.4}
\draw[white,thick] (0,1.4) to (1.8,-0.4);\draw[white,thick] (0,-1.4) to (1.8,0.4);
\draw[white,thick] (-\xx,\xx) to (\xx,-\xx);\draw[white,thick] (\xx,\xx) to (-\xx,-\xx);
\draw[draw=blue, very thick, densely dashed] (0,0) circle [radius=1cm];
\draw[red,thick] (-\xx,0) to (\xx,0);
\draw[black,fill=green] (-1,0) circle [radius=1.5mm];
\draw[black,fill=green] ( 1,0) circle [radius=1.5mm];
\end{tikzpicture}
\\[-2mm]
{\bf H1} & {\bf H2} & {\bf H3} & {\bf H4} & {\bf H5} & {\bf H6}
\end{tabular}
\end{center}
Thus $|Q\cap U|$
is equal to
$\infty$,
$4$, 
$2$, 
$4$, 
$2$ and 
$2$ at {\bf H1--H6}, \resp.

\newpage
We make an additional case distinction on the singular locus and real structure of 
a quadric surface $Q$:
\begin{Mlist}
\item[] {\bf S1.} $Q$ contains two real lines through each point.
\item[] {\bf S2.} $Q$ contains two complex conjugate lines through each point.
\item[] {\bf S3.}  $Q$ is singular.
\end{Mlist}
It follows from \EQN{q} that up to Euclidean similarity $q(0:x_1:x_2:x_3)$ is either 
\[
a_1\, x_1^2+a_2\, x_2^2+x_3^2,\quad 
a_1\, x_1^2+a_2\, x_2^2-x_3^2,\quad
a_1\, x_1^2+a_2\, x_2^2,\quad
a_1\, x_1^2-a_2\, x_2^2\quad\text{or}\quad 
x_1^2,
\]
for $a_1,a_2 \in \R_{>0}$. We separated the rows of \TAB{quad} according to these 5 cases.
If $a_i\neq 0$ in \EQN{q}, then $b_i=0$ after the translation $x_i\mapsto x_i-\frac{b_i}{2a_i}\,x_0$ for $1\leq i \leq 3$. 
After an Euclidean similarity we may assume that $b_0\in\{-1,0,1\}$.
If $a_i=0$ and $b_i\neq 0$ for $i\in\{2,3\}$, then $b_0=0$ after the translation $x_i\mapsto x_i-\frac{1}{b_0}\,x_0$. 
Notice that $a_1\neq a_2$, $a_1=a_2\neq 1$ and  $a_1=a_2=1$ are, 
a priori, treated as distinct cases \wrt {\bf H1--H5}.
It follows that the first column of \TAB{quad} lists, up to Euclidean similarity, all possible $q(x)$.

Complex circles in $\S^3$ are via a stereographic projection in 
one to one correspondence with complex lines and complex circles in $\P^3$. 
It follows that $|G(X)|$ corresponds to the number of 
pencils of complex lines or complex circles that cover $Q$.

Let us start with the row in \TAB{quad} corresponding to EH1. 
We see from the equation that we are in case {\bf H2} with {\bf S1} so that \#lines=2.
A pencil of planes, that contain a line spanned by two complex points in $Q\cap U$, 
defines a pencil of complex circles on $Q$. 
Thus $Q$ is covered by $\binom{4}{2}=6$ pencils of complex circles
so that $|G(X)|=6+2=8$.
There are two pairs of complex conjugate points in $Q\cap U$ and thus \#circles=2.

Let us now consider the row in \TAB{quad} corresponding to CH2.   
We see from the equation that we are in case {\bf H3} with {\bf S2} so that \#lines=0.
The pencil of planes, that contain the real line spanned by the two complex conjugate points in $Q\cap U$,
or the two complex conjugate lines that are tangent to both $Q$ and $U$,
define a pencil of complex circles that covers $Q$. 
It follows that \#circles=1 and $|G(X)|=2+3=5$.

\newpage
We treat the remaining rows analogously so that we obtain
the following table,
where M$i$:=empty-$i$ and  N$i$:=point-$i$:
\begin{center}
\footnotesize
\begin{tabular}{@{}c|c@{\hspace{5mm}}c@{\hspace{5mm}}c@{\hspace{5mm}}c@{\hspace{5mm}}c@{\hspace{5mm}}c@{}}
         & {\bf H1}  & {\bf H2}  & {\bf H3}  & {\bf H4} & {\bf H5} & {\bf H6}  \\\hline
{\bf S1} &           & EH1       & CH1       & HP       &          &           \\
{\bf S2} & $\S^2$,M3 & EE,M1,EH2 & CE,M2,CH2 & EP       & CP       &           \\
{\bf S3} & N3        & N1,EO     & N2,CO     & EY,N4,BY & CY,N5    &  PY       \\
\end{tabular}
\end{center}
The details of the remaining cases are now straightforward.
\end{proof}

\begin{lemma}
\label{lem:name}
If $X\subset\S^3$ is a celestial Darboux cyclide
with a real isolated singularity, then $\sigma_*$, $B(X)$ 
and its name are up to $\aut N(X)$ characterized by a row in \TAB{name}.
\end{lemma}

\begin{table}[!ht]
\caption{See \LEM{name}.}
\label{tab:name}
\centering
\begin{tabular}{lll|lll}
$\sigma_*$ & $B(X)$                        & name   & $\sigma_*$ & $B(X)$                               & name    \\\hline
$3A_1$     & $\{b_0\}$                     & EE/EH2 & $2A_1$     & $\{b_1,b_2,b_{12}\}$                 & EY    \\
$3A_1$     & $\{b_{13},b_{24}'\}$          & EP     & $2A_1$     & $\{b_1,b_2,b_{12},b_{13}',b_{24}'\}$ & CY    \\
$2A_1$     & $\{b_{12}\}$                  & EH1    & $2A_1$     & $\{b_{12},b_{34}\}$                  & EO    \\
$2A_1$     & $\{b_{13},b_{24},b_{12}'\}$   & CH1    & $2A_1$     & $\{b_{12},b_{34},b_{13}',b_{24}'\}$  & CO    \\  
$2A_1$     & $\{b_{12},b_{34}'\}$          & HP     &            &                                      & \\
\end{tabular}           
\end{table}

\begin{proof}
If the center of stereographic projection $\pi$ lies on the 
real isolated singularity of~$X(\R)$, then $Z:=\pi(X(\R))\subset \R^3$ is a quadric surface. 
We match $|G(X)|$ in \LEM{A3D4} and \LEM{D} with $|G(X)|$ in \PRP{Q}.
Moreover, if $W\subset B(X)$ \st $\sigma_*(W)=W$,
then it follows from \PRP{C}
that \#circles=$\set{c\in G(X)}{ \sigma_*(c)=c,~c\cdot W\nsucc 0}$
and
\#lines=$\set{c\in G(X)}{ \sigma_*(c)=c,~c\cdot W\succ 0}$.
For example, if $\sigma_*$  is $3A_1$ and $B(X)=\{b_0\}$, then
$G(X)=\{g_0,g_1,\g_{12},\g_{34},g_{13},g_{24},g_{14},g_{23}\}$
so that (\#circles, \#lines, $|G(X)|)=(2,0,8)$. 
Therefore, $Z$ is by \PRP{Q} either EE or EH2
so that $X$ is either an EE cyclide or an EH2 cyclide.
The remaining cases are similar so we conclude the proof of this lemma.
\end{proof}

\begin{remark}[Darboux cyclides]
\label{rmk:dp4}
Darboux cyclides \st $\sigma_*$ is $2A_1$ can be 
constructed with \ALG{construct}.
See \cite{tak2} for a systematic overview of implicit equations for Darboux cyclides
up to M\"obius transformations.
\newpage
If $(c_1,c_2,c_3,c_4)$ equals $(8,2,2,2)$ or $(8,2,2,-2)$, then
\[
Z:=\set{x\in\P^3}{(x_1^2+x_2^2+x_3^2)^2-x_0^2(c_1x_1^2-c_2x_2^2-c_3x_3^2-c_4x_0^2)=0}, 
\]
is a celestial surface of type $(2,4,3)$ \st $\sigma_*$ is $3A_1$ and~$D_4$, \resp, up to $\aut N(X)$.
See \citep[Figure~11]{pot2} 
and
\cite{nls-darboux} for Euclidean and kinematic constructions of Darboux cyclides, \resp.
\END
\end{remark}

\section{Combining the results}
\label{sec:results}

We prove the assertions in \SEC{intro} by applying the established results.

\begin{proof}[Proof of \THM{T}.]
If a Darboux cyclide is $\lambda$-circled, then $\lambda\leq 6$ by \LEM{sigma}.
Thus the listed types now follow from \PRP{dp} and \LEM{notype}.
It follows from \PRP{type} and \RMK{dp4}
that all listed types are realized by celestial surfaces.
\end{proof}

\COR{3} is a direct consequences of \THM{T}, 
but has already been established in \SEC{type}.

\begin{proof}[Proof of \COR{blum}.]
Direct consequence of \THM{T}.
\end{proof}

\begin{proof}[Proof of \THM{smooth}.]
Direct consequence of \PRP{dp} and \LEM{sigma}.
\end{proof}

\begin{proof}[Proof of \THM{B}.]
Direct consequence of the following assertions: \PRP{dp}, \LEM{sigma}, \LEM{A3D4}, \LEM{D} 
and \LEM{name}.
\end{proof}

\begin{proof}[Proof of \COR{graph}.]
Direct consequence of \THM{B} and \PRP{C}.
\end{proof}

\begin{proof}[Proof of \COR{line}.]
If $X\subset\S^n$ for $n\geq 3$ is a celestial surface that
contains infinitely many circles 
through some $p\in\S^n$, 
then either 
$X$ is a Veronese surface in $\S^4$, or
$p$ is a real base point of a pencil of circles.
In the latter case, the circle graph of $X$ contains a vertex
that is labeled with $+$ or $\times$.
If we choose the center of stereographic projection $\pi\c S^n\dto \R^n$ 
to be $p$,
then $\pi(X(\R))$ is a surface that is covered by circles and lines.
This corollary is thus a direct consequence of \THM{B} and \COR{graph}.
\end{proof}

\begin{proof}[Proof of \COR{sng}.]
We apply
\THM{B}, \LEM{A3D4}, \LEM{D}, \LEM{name}. 
Moreover,
$|E(X)|$ follows from \LEM{BGE},
since $X\cong X_N$ by \PRP{dp}.
\end{proof}

\begin{proof}[Proof of \COR{BEG}.]
Direct consequence of \THM{B} and \LEM{BGE}. 
\end{proof}

\begin{remark}
\label{rmk:dc}
The \df{cyclicity} of a surface $Z\subset\R^3$ is defined as the multiplicity of 
the absolute conic~$U$
in $\P(Z)$ (recall \DEF{P} and \DEF{U}). 
Notice that if $\P(Z)=\set{x\in\P^3}{ F(x)=0}$ is of cyclicity~$c$
and degree~$d$, 
then $F(0:x_1:x_2:x_3)=(x_1^2+x_2^2+x_3^2)^c\,L(x)$
for some form~$L(x)$ of degree~$d-2\,c$. 
It now follows from \PRP{dc} below that our definition for ``Darboux cyclide'' 
is equivalent to the definition in \citep[Section~2]{pot2}.
\END
\end{remark}

\begin{proposition}
\label{prp:dc}
If $Z\subset\R^3$ is a surface of degree~$d$ and cyclicity~$c$, 
then 
$\deg \pi^{-1}(Z)=2\,(d-c)$, 
for all stereographic projections $\pi\c S^3\dto \R^3$.
\end{proposition}

\begin{proof}
We may assume up to $\aut(S^3)$ that $\pi$ has center $(0,0,0,1)$.
Let $\piP\c\S^3\dto\P^3$ be its projective closure 
so that
$\piP(x)=(x_0-x_4:x_1:x_2:x_3)$ 
and  $\piP^{-1}(y)=(\Delta+y_0^2:2y_0y_1:2y_0y_2:2y_0y_3:\Delta-y_0^2)$
with $\Delta:=y_1^2+y_2^2+y_3^2$.
Let $X\subset \S^3$ denote the M\"obius model of $Z$
so that 
$\deg X=\deg \pi^{-1}(Z)$
and so that
$\P(Z)$ coincides with the Zariski closure of~$\piP(X)$.
We consider the linear series associated to $\piP^{-1}$, namely the quadric surfaces in~$\P^3$
that are preimages of hyperplanes in $\P^4$ parametrized by $\alpha\in\P^4$:
\[
Q_\alpha:=\set{y\in \P^3}{
\alpha_0\,(\Delta+y_0^2)+  
\alpha_1\,y_0\,y_1+
\alpha_2\,y_0\,y_2+
\alpha_3\,y_0\,y_3+
\alpha_4\,(\Delta-y_0^2)=0}.
\]
The degree of $X\subset\P^4$ is by definition the number of common intersections
between $X$ and two general hyperplanes \citep[Section~I.7]{har}.
Notice that $Q_\alpha$ contains the absolute conic~$U$ for all $\alpha\in \P^4$ 
and thus $\piP^{-1}$ is not defined at $U$.
Hence,  
the degree of $X$ is equal to the number of complex points in
$Q_{\alpha_0}\cap Q_{\alpha_1}\cap X$ that do not lie in $U$,
for some general $\alpha_0,\alpha_1\in\P^4$.
The complex conic $C$ defined by $Q_{\alpha_0}\,\cap\,Q_{\alpha_1}$
intersects $X$ at $U$ in $2\,c$ points when counted 
with multiplicity (see \citep[Section~2.1]{intersection}).
Since $C$ intersects $X$ in $2\,d$ points by B\'ezout's theorem
it follows that $\deg X=2\,d-2\,c$ as was to be shown.
We remark that if $|C(\R)|=\infty$, then $\iota^{-1}(C(\R))\subset\R^3$ is a circle (see \DEF{P})
and thus $2\,(d-c)$ is equal to 
the maximal number of intersections of $Z\subset\R^3$ with a circle.
\end{proof}

\begin{definition}
\label{def:cweb}
A \df{complex hexagonal web} is defined 
by replacing
the points, curves and surfaces in \DEF{web}
with complex points, complex curves and complex surfaces, \resp.
\END
\end{definition}

\newpage
\begin{example}
\label{exm:L}
Suppose that 
$q_1,q_2,q_3\in \P^2$ are distinct complex points that do not lie on a complex line.
It follows from \cite{graf1} 
that the following set forms a complex hexagonal web of lines on $\P^2$ (see \DEF{cweb}):
\[
\cLL:=\set{L\subset \P^2}{L \text{ is a complex line \st }|L\cap\{q_1,q_2,q_3\}|\geq 1}.
\]
See alternatively \citep[Section~1.3, page~24]{bla1}.
\END
\end{example}

\begin{proof}[Proof of \THM{hex}.]
Suppose that $X\subset\S^n$ is the M\"obius model of~$Z\subset\R^n$
so that $Z=\pi(X(\R))$ is $\lambda$-circled with $\lambda\geq 3$.
It follows from \COR{graph} that $Z$ is not covered by lines,
because the graphs in \TAB{graph} 
containing an $\oplus$ vertex contain less than 3 vertices 
that are not labeled by $\oplus$.
Hence, $Z$ is covered
by hexagonal web of circles if and only if $X$ is 
covered by a hexagonal web of circles.

Let $\cLL$ be the complex hexagonal web of lines in $\P^2$ from \EXM{L}
that pass through $q_1$, $q_2$ or $q_3$.
Our plan is to construct a complex
birational map~$\mu\c X\dto \P^2$ and a 3-web~$\cW$ of
circles on~$X$ \st $\mu(\cW)\subset\cLL$ and $|C\cap C'|\leq 1$ for all $C,C'\in\cW$.
In particular, we require that $\mu$ is almost everywhere a complex isomorphism
that sends circles in~$\cW$ to complex lines in~$\cLL$.
Notice that circles in $\cW$ are either disjoint or intersect in a single point 
that must be real.
We follow the procedure in \RMK{web} and when we draw distinct circles $C,C'\in \cW$
we also draw complex lines $\mu(C),\mu(C')\in\cLL$.
If $\mu(C)\cap \mu(C')\cap \{q_1,q_2,q_3\}=\emptyset$, 
then $C$ and $C'$ must intersect in a point~$p$ and thus there exists a circle $C''\in\cW$
so that $p\in C''$ and $C''\notin\{C,C'\}$. 
Hence the procedure results in a closed hexagon so that
$\cW$ must by \DEF{web} be a hexagonal web of circles.

We make a case distinction on the 
possible types of~$X$ as listed in \THM{T}.

First suppose that $X$ is of type either $(3,6,4)$ or $(3,6,5)$. 
We let $\varphi\c Y\to X$ denote its smooth model.
We know from \PRP{C}a that $|G(X)|\geq 3$.
Recall from \LEM{BGE}[a,c] and \PRP{dp} that 
$N(X)\cong\bas{\l_0,\l_1,\p_1,\p_2}_\Z$,
$G(X)=\{g_0,g_1,g_{12}\}$ and $B(X)=\emptyset$.
It follows that~$Y$ is isomorphic to the blowup of~$\P^1\times\P^1$ 
in two complex conjugate centers~$p_1$ and~$p_2$ (see \RMK{S1S1}). 
Thus $Y$ is also complex isomorphic to the 
blowup of $\P^2$ in three complex centers that are not collinear.
We may assume up to $\aut_\C\P^2$ that these centers coincide with 
$q_1$, $q_2$ and $q_3$ in \EXM{L}.
Let $\tau\c Y\to\P^2$ be the complex birational morphism that 
contracts three complex $(-1)$-curves to these centers.
We know from 
\citep[Example~II.8.20.3, Proposition~V.3.2, Proposition~V.3.3]{har}
that
$N(X)\cong\bas{\alpha,\beta_1,\beta_2,\beta_3}_\Z$
and $-\k=3\,\alpha-\beta_1-\beta_2-\beta_3$ is the anticanonical class of~$Y$.
Notice that $\alpha$ is the pullback via~$\tau$ of the class of 
a general complex line in~$\P^2$ and that $\beta_i$
is the class of the complex $(-1)$-curve that is centered at $q_i$ for $1\leq i\leq 3$.
Without loss of generality we have 
$\alpha=\l_0+\l_1-\p_1$, $\beta_1=\l_1-\p_1$, $\beta_2=\l_0-\p_1$
and $\beta_3=\p_2$.
It follows that $\alpha-\beta_i$ is the class of a curve in $Y$ that 
is send by~$\tau$ to a complex line in~$\P^2$ that passes through~$q_i$
(alternatively see \citep[Proposition~V.3.6]{har}).
We remark that the image of the complex $(-1)$-curve in $Y$ with class $\alpha-\beta_i-\beta_j$
is the complex line in $\P^2$ that passes through $q_i$ and $q_j$
for all $1\leq i<j\leq 3$.
Now let $\mu$ be defined as the complex birational map $\tau\circ\varphi^{-1}\c X\dto \P^2$ 
and let $\cW$ be the set of circles in~$X$.
Since
$G(X)=\{g_0,g_1,g_{12}\}=\{\alpha-\beta_1,\alpha-\beta_2,\alpha-\beta_3\}$,
it follows that $\mu(\cW)\subset\set{L\in \cLL}{ |L\cap\{q_1,q_2,q_3\}|=1}$
and $|C\cap C'|=[C]\cdot [C']\leq 1$ for all $C,C'\in\cW$. 
Hence, $X$ is covered by a hexagonal web of circles
so that $Z$ is covered by a hexagonal web of circles as well.

Next suppose that $X$ is of type~$(\infty,4,4)$
so that $X\subset \P^5$ is a Veronese surface by \THM{smooth}.
By definition there exists a biregular isomorphism~$\mu\c X\to \P^2$
such that circles in $X$ are send to lines in $\P^2$.
We may assume up to $\aut_\C\P^2$ that $\mu^{-1}(\{q_1,q_2,q_3\})$
consists of three points which are real.
We define $\cW$ to be the set of circles in $X$
\st $\mu(\cW)\subset\cLL$ and thus $|C\cap C'|=1$ for all $C,C'\in\cW$.
This concludes the proof for this case as $Z$ must 
be covered by a hexagonal web of circles.

Now suppose that $X$ is of type~$(\infty,2,2)$ so that $X$ is a sphere.
In this case the inverse stereographic projection of a hexagonal 
web of lines in the plane is a hexagonal web of circles in~$X$.

Finally suppose that $X$ is of type $(\lambda,4,3)$
so that $Z$ is a Darboux cyclide.
Recall from \RMK{dc} that
our definition of Darboux cyclide 
coincides with the definition in \cite{pot2}.
Thus in this case it follows from 
\citep[Theorem~18]{pot2} together with
\RMK{web2} and \COR{graph}
that $Z$ is covered by a hexagonal web of circles. 

We concluded the proof, since by \THM{T}
we considered all possible types $(\lambda,n,d)$ of $X$ \st $\lambda\geq 3$.
\end{proof}

\section{Acknowledgements}

I would like to thank R.~Krasauskas, H.~Pottmann, J.~Schicho, M.~Skopenkov and S.~Zub\.{e}
for the inspiring discussions, which have have been invaluable for this paper.
In particular, the detailed comments of M.~Skopenkov were extremely helpful.
I thank J.~Koll\'ar for interesting historical remarks.
The computations were done using \citep[Sage]{sage} and \citep[Magma]{magma}.
The images were made using \citep[Povray]{povray}.
This research was supported 
by base funding of the King Abdullah University of Science and Technology (KAUST)
and by the Austrian Science Fund (FWF) project P33003.

\bibliography{fam-circles}

\paragraph{address of author:}
Johann Radon Institute for Computational and Applied Mathematics (RICAM), 
Austrian Academy of Sciences
\\
\textbf{email:} niels.lubbes@gmail.com

\end{document}